\documentclass[preprint,12pt]{elsarticle}

\usepackage{amssymb}
\usepackage{amsmath}
\usepackage{float}
\usepackage{pgfplots}
\pgfplotsset{compat=newest}

\journal{Journal of Computational Physics}

\begin{document}

\begin{frontmatter}

\title{AP-Cloud: Adaptive Particle-in-Cloud Method for Optimal Solutions to Vlasov-Poisson Equation}
\date{\today} 

\author[label1]{Xingyu Wang}
\author[label1,label2]{Roman Samulyak\footnote{Corresponding author, roman.samulyak@stonybrook.edu}}
\author[label1]{Xiangmin Jiao}
\author[label1,label2]{Kwangmin Yu}

\address[label1]{Department of Applied Mathematics and Statistics,
  Stony Brook University, Stony Brook, NY 11794}
\address[label2]{Computational Science Center, Brookhaven National
  Laboratory, Upton, NY 11973}

\begin{abstract}
We propose a new adaptive Particle-in-Cloud (AP-Cloud) method for obtaining optimal numerical solutions to the Vlasov-Poisson equation. Unlike the traditional particle-in-cell (PIC) method, 
which is commonly used for solving this problem, the AP-Cloud adaptively selects computational
nodes or particles to deliver higher accuracy and efficiency when the particle distribution is highly non-uniform. Unlike other adaptive techniques for PIC, our method balances the errors 
in PDE discretization and Monte Carlo integration, and discretizes the differential 
operators using a generalized finite difference (GFD) method based on a weighted least 
square formulation. As a result, AP-Cloud is independent of the geometric shapes of 
computational domains and is free of artificial parameters. Efficient and robust 
implementation is achieved through an octree data structure with 2:1 balance. We analyze 
the accuracy and convergence order of AP-Cloud theoretically, and verify the method using 
an electrostatic problem of a particle beam with halo. Simulation results show that 
the AP-Cloud method is substantially more accurate and faster than the traditional PIC, 
and it is free of artificial forces that are typical for some adaptive PIC techniques.

\end{abstract}

\begin{keyword}
particle method \sep generalized finite difference \sep PIC \sep AMR-PIC

\MSC 65M06 \sep 70F99 \sep 76T10

\end{keyword}
\end{frontmatter}

\section{Introduction}

The Particle-in-Cell (PIC) method \cite{HockneyEastwod} is a popular method for solving the Vlasov-Poisson equations for a class of problems 
in plasma physics, astrophysics, and particle accelerators, for which electrostatic approximation applies, as well as for solving the gravitational problem in 
cosmology and astrophysics.  In such a hybrid particle-mesh method, the distribution function is approximated using particles and the Poisson problem is 
solved on a rectangular mesh. Charges (or masses) of particles are interpolated onto the mesh, and the Poisson problem is discretized using finite differences or
spectral approximations. On simple rectangular domains, FFT methods are most commonly used for solving the Poisson problem. In the presence of irregular 
boundaries, finite difference approximations are often used, complemented by a cut-cell (a.k.a. embedded boundary) method \cite{Colella_EBM} for computational cells 
near boundaries, and fast linear solvers (including multigrid iterations) for the corresponding linear system. The computed force (gradient of the potential)
on the mesh is then interpolated back to the location of particles. For problems with irregular geometry, unstructured grid with finite element method is often used.

The traditional PIC method has several limitations. It is less straightforward to use for geometrically complex domains. The aforementioned embedded boundary 
method, while maintaining globally second order accuracy for the second order finite difference approximation, usually results in much larger errors 
near irregular boundaries \cite{WangSam10}. It is also difficult to generalize to higher order accuracy. 

Another major drawback of the PIC method is associated with highly non-uniform distribution of particles. As shown in Section 2, the discretization of the 
differential operator and the right hand side in the PIC method is not balanced in terms of errors. The accuracy is especially degraded in the presence of non-uniform 
particle distributions. The AMR-PIC \cite{VayCol1,VayCol2} improves this problem by performing block-structured adaptive mesh refinement of a rectangular mesh, so 
that the number of particles per computational cell is approximately the same. However, 
the original AMR-PIC algorithms suffered from very strong artificial self-forces due to spurious
images of particles across boundaries between coarse and refined mesh patches. Analysis of self-force sources and a method for their mitigation was proposed in
\cite{ColNor10}. 

In this paper, we propose a new  adaptive Particle-in-Cloud (AP-Cloud) method for obtaining optimal numerical solutions to the
Vlasov-Poisson equation. Instead of a Cartesian grid as used in the traditional PIC, the
AC-Cloud uses adaptive computational nodes or particles with an octree data structure.
 The quantity characterizing particles (charge in electrostatic problems or 
 mass in gravitational problems) is assigned to computational nodes by a weighted least squares approximation. The partial differential 
 equation is then discretized using a generalized finite difference (GFD) method and solved with fast linear solvers.  The density of 
 nodes is chosen adaptively, so that the error from GFD and that from Monte Carlo integration are balanced, and the 
 total error is approximately minimized. The method is independent of geometric shape of computational domains and free of artificial self-forces.

The remainder of the paper is organized as follows. In Section 2, we analyze numerical errors of the traditional PIC method and formulate optimal refinement strategy.
The AP-Cloud method, generalized finite differences, and the relevant error analysis are presented in Section 3. Section 4 describes some implementation details of the method. 
Section 5 presents numerical verification tests using 2D and 3D problems of particle beams with halo and additional tests demonstrating the absence of artificial self-forces. 
We conclude this paper with a summary of our results and perspectives for the future work.

\section{Error analysis of particle-in-cell method}
In a particle-in-cell (PIC) method, the computational objects include a large number of 
particles and an associated Cartesian grid. These particles are typically 
randomly sampled and represent an even greater number of physical particles
(e.g., protons), so they are also known as ``macro-particles'', but conventionally
simply referred to as particles. For simplicity of presentation, we will focus on 
electrostatic problems, for which the states are particle charges. Suppose there are 
$N$ charged particles at positions $\mathcal{P}=\{\mathbf{p}^i\mid i=1,2,\dots,N\}$ in $D$ dimensions, and let $q^i$ denote the charge 
at $\mathbf{p}^i$. For simplicity, we assume that all the particles carry the same amount of charge, 
the total charge is 1, i.e., $Nq^i=\int_\Omega\rho(\mathbf{x})\mathbf{dx}=1$, and
the charges of the particles can be represented accurately by a continuous 
charge distribution function $\rho$. We assume that $\rho$ is smooth and positive, 
and its value and all derivatives have comparable magnitude.
Let $M$ denote the Cartesian grid, and without loss of generality, suppose its
edge length is $h$ along all directions, and let
$\mathbf{y}^j$ denote the $j$th grid point in $M$.
A PIC method estimates the charge density $\rho$ on $M$, then solves the Poisson equation 
\begin{equation}
\Delta \phi=c\rho
\label{eq:Poisson}
\end{equation}
on $M$ to obtain the potential $\phi$, whose gradient is the electric field $\mathbf{E}$.
In this setting, a PIC method consists of the following three steps:
\begin{enumerate}
\item Approximate the right-hand side of (\ref{eq:Poisson}) by interpolating the states 
from particles to the grid points $\mathbf{y}^j$, i.e.,
\begin{equation}\label{eq:density}
\tilde\rho(\mathbf{y}^j,\mathcal{P},h)=\frac{1}{h^D}
\sum_{i=1}^N q^i\Phi\left(\frac{\mathbf{p}^i-\mathbf{y}^j}{h}\right)
\approx \underbrace{\frac{1}{h^{D}}\int_{\Omega}\rho(\mathbf{x})\Phi\left(\frac{\mathbf{x}-\mathbf{y}^j}{h}\right)\mathbf{dx}}_{\bar{\rho}(\mathbf{y}^{j},h)}
\approx\rho(\mathbf{y}^j),
\end{equation}
where $\Phi$ is the interpolation kernel, a.k.a. the charge assignment scheme. 
\item Discretize the left-hand side of (\ref{eq:Poisson}) on $M$, 
typically using the finite difference method, and then solve
the resulting linear system. 
\item Obtain the electric field $\mathbf{E}$ by computing $\nabla\phi$ using 
finite difference,
and then interpolating $\mathbf{E}$ from the grid points to the particles using the same
interpolation kernel $\Phi$ as in Step 1, i.e.,
\begin{equation}
\mathbf{\tilde E}(\mathbf{p}^i,\mathcal{P},h) =\frac{1}{h^D}\sum_{\mathbf{y}^j\in M}\mathbf{\tilde E}(\mathbf{y}^j,\mathcal{P},h)\Phi\left(\frac{\mathbf{p}^i-\mathbf{y}^j}{h}\right)\approx\mathbf{E}(\mathbf{p}^i).
\end{equation}
\end{enumerate}
One of the most commonly used charge assignment schemes $\Phi$ is the 
cloud-in-cell (CIC) scheme
\begin{equation}\label{cloud_in_cell}
\Phi(\mathbf{x})=\prod_{d=1}^D\max{(1-|x_d|,0)},
\end{equation}
for which the interpolation in Step 3 corresponds to bilinear and trilinear interpolation
in 2-D and 3-D, respectively.

In PIC, the error in potential $\phi$ comes from two sources. One is from the first approximation in (\ref{eq:density}),
for which the analysis is similar to Monte Carlo integration within a control volume 
associated with $\mathbf{y}^j$, under the assumption that (\ref{eq:density})
is a continuous function. The other source is the discretization error of both the second approximation in (\ref{eq:density}) in Step 1 and the left-hand
side of (\ref{eq:Poisson}) on $M$ in Step 2. 
We denote the above two errors from these two sources as 
$\mathcal{E}_{M}$ and $\mathcal{E}_{D}$, respectively. As shown in \ref{appen}, under the assumption that the
interpolation kernel $\Phi$ satisfies the positivity condition, 
the expected value of the former is 
$$\mbox{E}[|\mathcal{E}_{M}|]=\mathcal{O}\left(\sqrt{\frac{\rho(\mathbf{y})}{Nh^D}}\right),$$
and the discretization
 error is
$$\mathcal{E}_{D}=\mathcal{O}\left(\rho(\mathbf{y})h^2\right).$$
Let $A$ denote the coefficient matrix of the linear system in step 2, and suppose 
$\Vert A^{-1}\Vert$ is bounded by a constant. The expected total error in the 
computed potential $\phi$ is then 
$$\mathcal{E}=\mathcal{O}(\mbox{E}[|\mathcal{E}_{M}|])+\mathcal{O}(\mathcal{E}_{D})=
\mathcal{O}\left(\sqrt{\frac{\rho(\mathbf{y})}{Nh^D}}+\rho(\mathbf{y})h^2\right).$$

In general, $\mathcal{E}_{D}$ dominates the total error for coarse grids and $\mathcal{E}_{M}$ dominates 
for finer grids. The total expected error is approximately minimized if 
$\mathcal{E}_{M}$ and $\mathcal{E}_{D}$ are balanced. If the particles are uniformly distributed, 
then the errors are balanced when 
\begin{equation}\label{eq:optimal}
h=\mathcal{O}\left({N\rho(\mathbf{y})}\right)^{-\frac{1}{4+D}}.
\end{equation}
In this setting, the discretization error in $\phi$ is second order in $h$. The discretization error in numerical differentiation $\mathbf{E}=\nabla\phi$ is also second order, a fact called supraconvergence~\cite{sc}. Thus, although the optimal mesh size is deduced to minimize the error in $\phi$, the error in $\mathbf{E}$ is also minimized. 

In many applications, the particle distribution is highly non-uniform, for which 
the PIC is neither efficient nor accurate. In \cite{VayCol2}, 
an adaptive method, called AMR-PIC, was
proposed, which fixed the number of particles per cell, and hence
$$h=\mathcal{O}\left(N\rho(\mathbf{y})\right)^{-\frac{1}{D}}.$$
The AMR-PIC over-refines the grid compared to the optimal grid resolution in (\ref{eq:optimal}). 
In addition, the original AMR-PIC technique also introduces artificial self forces. New adaptive
strategies are needed to resolve both of these issues.

\section{Adaptive Particle-in-Cloud method}
In this section, we describe a new adaptive method, called Adaptive Particle-in-Cloud or 
AP-Cloud, which approximately minimizes the error by balancing Monte Carlo noise and discretization error, and at the same time is free of the artificial self forces present 
in AMR-PIC.

The AP-Cloud method can be viewed as an adaptive version of PIC that replaces the traditional Cartesian mesh of PIC by an octree data
structure. We use a set of computational nodes, which are octree cell centres, instead of the Cartesian grid, of which the distribution is derived using an error balance criterion. Computational nodes will be referred to as nodes in the remainder of the paper. Instead of the finite difference discretization of the Laplace operator, we use the method of generalized finite-difference (GFD) \cite{Benito03}, based on a weighted least squares formulation. 
The framework includes interpolation, least squares approximation, and numerical differentiation on a stencil in the form of cloud of nodes in a neighborhood 
of the point of interest. It is used for the charge assignment scheme, numerical differentiation, and interpolation of solutions. 
The advantage of GFD is that it can treat coarse regions, refined regions, and refinement boundaries in the same manner, and it is more flexible for problems in complex domain or with irregular refinement area. As a method of integration, GFD will be used in the quadrature rule in charge assignment scheme.
The new charge assignment scheme, together with GFD differentiation and interpolation operators from computational nodes to particles is easily generalizable to higher 
order schemes. We described the key components of the AP-Cloud in this section, and
discuss its implementation details in Section~\ref{sec:implementation}.

\subsection{Generalized finite-difference method}

For simplicity of presentation, we consider a second order generalized finite-difference method.

Let $\mathbf{y}^j$, $j\in\{1,2,\cdots,m\}$ be the nodes in a neighborhood of reference node $\mathbf{y}^0$. Given a $C^2$ function $f$, by Taylor expansion we have
\begin{equation}
f(\mathbf{y}^j)=f(\mathbf{y}^0)+(\mathbf{y}^j-\mathbf{y}^0)^T\nabla f(\mathbf{y}^0)+(\mathbf{y}^j-\mathbf{y}^0)^T H(\mathbf{y}^0)(\mathbf{y}^j-\mathbf{y}^0)+O(h^3),
\end{equation}
where $h$ is the characteristic interparticle distance in the neighborhood, for example, $h=\max\limits_{d,j}|y_d^j-y_d^0|$, and $H$ is the Hessian matrix.
Putting equations for all neighbors together and omitting higher order term, we obtain
\begin{equation}\label{gfd}
V(\mathbf{y}^0)\partial f(\mathbf{y}^0)=\delta f(\mathbf{y}^0),
\end{equation}
where $V(\mathbf{y}^0)$ is a generalized Vandermonde matrix, $\partial f(\mathbf{y}^0)$ is the first order and second order derivative of $f$ at $\mathbf{y}^0$, and $\delta f(\mathbf{y}^0)$ is the increment of $f$.

For example, in 2-D, let $m=5$ and $\delta y_d^j=y_d^j-y_d^0$, then
\begin{equation}
V(\mathbf{y}^0)=
\left[ 
\begin{array}{ccccc}
\delta y_1^1 & \delta y_2^1 & \frac{1}{2}(\delta y_1^1)^2 & \delta y_1^1 \delta y_2^1 & \frac{1}{2}(\delta y_2^1)^2\\
\delta y_1^2 & \delta y_2^2 & \frac{1}{2}(\delta y_1^2)^2 & \delta y_1^2 \delta y_2^2 & \frac{1}{2}(\delta y_2^2)^2\\
\delta y_1^3 & \delta y_2^3 & \frac{1}{2}(\delta y_1^3)^2 & \delta y_1^3 \delta y_2^3 & \frac{1}{2}(\delta y_2^3)^2\\
\delta y_1^4 & \delta y_2^4 & \frac{1}{2}(\delta y_1^4)^2 & \delta y_1^4 \delta y_2^4 & \frac{1}{2}(\delta y_2^4)^2\\
\delta y_1^5 & \delta y_2^5 & \frac{1}{2}(\delta y_1^5)^2 & \delta y_1^5 \delta y_2^5 & \frac{1}{2}(\delta y_2^5)^2\\
\end{array}
\right],
\end{equation}
\begin{equation}
\partial f(\mathbf{y}^0)=\left[f_{y_1}(\mathbf{y}^0),f_{y_2}(\mathbf{y}^0),f_{y_1 y_1}(\mathbf{y}^0),f_{y_1 y_2}(\mathbf{y}^0),f_{y_2 y_2}(\mathbf{y}^0)\right]^T,
\end{equation}
where $f_{y_d}(\mathbf{y}^0)$ denotes the derivative of $f$ with respect to $y_d$, and
\begin{equation}
\delta f(\mathbf{y}^0)=\left[f(\mathbf{y}^1)-f(\mathbf{y}^0),f(\mathbf{y}^2)-f(\mathbf{y}^0),\dots,f(\mathbf{y}^5)-f(\mathbf{y}^0)\right]^T.
\end{equation}

To analyze the error in GFD, let $\xi_d^j=h^{-1}\delta y_d^j$. Rewrite (\ref{gfd}) as
\begin{equation}\label{regfd}
V_0(\mathbf{y}^0)\partial f_0(\mathbf{y}^0)=\delta f(\mathbf{y}^0),
\end{equation}
where 
\begin{equation}
V_0(\mathbf{y}^0)=
\left[ 
\begin{array}{ccccc}
\xi_1^1 & \xi_2^1 & \frac{1}{2}({\xi_1^1})^2 & \xi_1^1\xi_2^1 & \frac{1}{2}({\xi_2^1})^2 \\
\xi_1^2 & \xi_2^2 & \frac{1}{2}({\xi_1^2})^2 & \xi_1^2\xi_2^2 & \frac{1}{2}({\xi_2^2})^2 \\
\xi_1^3 & \xi_2^3 & \frac{1}{2}({\xi_1^3})^2 & \xi_1^3\xi_2^1 & \frac{1}{2}({\xi_2^3})^2 \\
\xi_1^4 & \xi_2^4 & \frac{1}{2}({\xi_1^4})^2 & \xi_1^4\xi_2^4 & \frac{1}{2}({\xi_2^4})^2 \\
\xi_1^5 & \xi_2^5 & \frac{1}{2}({\xi_1^5})^2 & \xi_1^5\xi_2^5 & \frac{1}{2}({\xi_2^5})^2 \\
\end{array}
\right],
\end{equation}
\begin{equation}
\partial f_0(\mathbf{y}^0)=\left[hf_{y_1}(\mathbf{y}^0),hf_{y_2}(\mathbf{y}^0),h^2f_{y_1 y_1}(\mathbf{y}^0),h^2f_{y_1 y_2}(\mathbf{y}^0),h^2f_{y_2 y_2}(\mathbf{y}^0)\right]^T.
\end{equation}
Now $V_0(\mathbf{y}^0)$ depends on the shape but not the diameter of the GFD stencil. The error in solving linear system is
\begin{equation}
\Vert\mbox{Error}(\partial f_0(\mathbf{y}^0))\Vert\leq \Vert V_0^{-1}(\mathbf{y}^0)\Vert \Vert \mbox{Error} (\delta f(\mathbf{y}^0))\Vert.
\end{equation}
The error in right hand size comes from the omitted term in the Taylor expansion $\mathcal{O}(f_{y_iy_jy_k}h^{3})=\mathcal{O}(fh^{3})$ for $i,j,k\in{1,2}$, and $\Vert V_0^{-1}(\mathbf{y}^0)\Vert$ is a constant independent of $h$, so the error in $\partial f_0(\mathbf{y}^0)$ is also $\mathcal{O}(fh^{3})$. Because the coefficient before the $l$th order derivative in $\partial f_0(\mathbf{y}^0)$ is $h^l$, the error for $l$th order derivative is $\mathcal{O}(fh^{3-l})$

In this example, the number of neighbors is equal to the number of unknowns, and it is quite
likely for $V_0(\mathbf{y}^0)$ to be nearly singular. In practice, we use more neighbors in the stencil than the number of coefficients in the Taylor series to improve the stability. In AP-Cloud method, 8 neighbors instead of 5 are used for the second order GFD in two dimensions, and 17 neighbors instead of 9 are used in three dimensions. In this case, the linear system is a least square problem. It is often helpful to assign more weights to closer neighbors to improve the accuracy, which is called the weighted lest square method~\cite{Benito01}. AP-Cloud method uses a normalized Gaussian weight function~\cite{Onate}: 
\begin{equation}
W(r)=\frac{e^{-r^2/r_{max}^2}-e^{-c}}{1-e^{-c}},
\end{equation}
where $W$ is the weight, $r$ is the distance of the neighbor from the reference particle, $r_{max}$ is the maximum distance of all neighbors in the stencil from the reference particle, $c=4$.

By solving the linear system or least square problem (\ref{gfd}), we can express the gradient $\partial f$ as linear combinations of $\delta f$. For example, once the potential $\phi$ is computed at nodes, we can find its gradient by generalized finite-difference, and then interpolate it to particles by Taylor expansion. Generally, the error of the $k$th order GFD interpolation is $\mathcal{O}(f(\mathbf{y})h^{k+1})$, and its approximation of the $l$th order derivative is $\mathcal{O}(f(\mathbf{y})h^{k-l+1})$. 

Given a set of nodes, the selection for GFD neighbors, or the shape of the GFD stencil, is important for both accuracy and stability. Simply choosing the nearest nodes to be neighbors may lead to an imbalanced stencil. We follow the quadrant criterion in~\cite{3} and select two nearest nodes from each quadrant to be neighbors.

\subsection{Algorithm of AP-Cloud method}
AP-Cloud also has three steps to calculate the electric field given by a particle distribution: a density estimator, a Poisson solver, and an 
interpolation step, but each step is different from its counterpart in PIC. Let $\mathcal{C}$ be the set of all computational nodes, and $f(\mathcal{C})=(f(\mathbf{y}^1),f(\mathbf{y}^2),\cdots,f(\mathbf{y}^n))$, where $n$ is the total number of nodes. Below is a detailed description
of the three steps.

\begin{enumerate}
\item Approximate density by interpolating states 
from particles $\mathcal{P}$ to computational nodes $\mathcal{C}$ by
\begin{eqnarray}\label{eq:density_apcloud}
&&a(\mathbf{y}^j)^T\partial \tilde{\rho}(\mathbf{y}^j)\\ \nonumber
:&=&\sum_{l=0}^k\sum_{d_1,d_2,\cdots,d_l=1}^D\underbrace{\frac{1}{h^D}\int_\Omega\prod_{i=1}^l(x_{d_i}-y_{d_i}^j)\Phi\left(\frac{\mathbf{x}-\mathbf{y}^j}{h}\right)\mathbf{dx}}_{a_{d_1,d_2,\cdots,d_l}(\mathbf{y}^j)}\tilde{\rho}_{y_{d_1}y_{d_2}\cdots y_{d_l}}(\mathbf{y}^j)\\ \nonumber
&=&\underbrace{\frac{1}{h^D}\sum_{i=1}^N q^i\Phi\left(\frac{\mathbf{p}^i-\mathbf{y}^j}{h}\right)}_{\rho_{M}(\mathbf{y}^j,h)}.
\end{eqnarray}

The right hand side  of (\ref{eq:density_apcloud}) is identical to the Monte Carlo integration in PIC, but the left hand side is a linear combination of derivatives, instead of the simple $\tilde\rho(\mathbf{y}^j,\mathcal{P},h)$ in the PIC method. Because the coefficients in the linear combination, $a(\mathbf{y}^j)=(a_{d_1,d_2,\cdots,d_l}(\mathbf{y}^j))$, depend only on $h$ and the interpolation kernel, they can be easily pre-calculated and tabulated in a lookup table. The derivatives $\tilde{\rho}_{y_{d_1}y_{d_2}\cdots y_{d_l}}(\mathbf{y}^j)$ are in turn linear combinations of density values $\rho(\mathbf{y}^j)$ given by least square solution of  (\ref{gfd}):
\begin{equation}\label{eq:pseudoinverse}
\partial \tilde{\rho}(\mathbf{y}^j)=V(\mathbf{y}^j)^{+}\delta \tilde{\rho}(\mathbf{y}^j),
\end{equation}
where $V(\mathbf{y}^j)^{+}$ is the pseudo-inverse of the Vandermonde matrix.

Let $C(\mathbf{y}^j)$ be a matrix such that $\delta\tilde{\rho}(\mathbf{y}^j)=C(\mathbf{y}^j)\tilde{\rho}(\mathcal{C})$. Substituting (\ref{eq:pseudoinverse}) into (\ref{eq:density_apcloud}), we get a linear equation for density values at the reference node $\mathbf{y}^j$ and its neighors
\begin{equation}
\underbrace{a(\mathbf{y}^j)^TV(\mathbf{y}^j)^{+}C(\mathbf{y}^j)}_{b(\mathbf{y}^j)}\tilde{\rho}(\mathcal{C})=\rho_{M}(\mathbf{y}^j,h).
\end{equation}

Putting equations for all nodes together, we obtain a global linear system for density values
\begin{equation}\label{eq:density_linear_system}
B(\mathcal{C})\tilde{\rho}(\mathcal{C})=\rho_{M}(\mathcal{C},h).
\end{equation}
where $B(\mathcal{C})=[b(\mathbf{y}^1);b(\mathbf{y}^2);\cdots;b(\mathbf{y}^n)]$. Solution of (\ref{eq:density_linear_system}) is the estimated density in AP-Cloud method.
\item Discretize the left-hand side of (\ref{eq:Poisson}) on $\mathcal{C}$ using GFD method. Solve the resulting linear system for $\phi(\mathcal{C})$.
\item Obtain the electric field $\mathbf{E}$ by computing $\nabla\phi$ using GFD method,
and then interpolating $\phi$ and $\mathbf{E}$ from the $\mathcal{C}$ to $\mathcal{P}$ using a Taylor expansion.
\end{enumerate}

\subsection{Error analysis for AP-Cloud}
Similar to PIC, the error in AP-Cloud also contains the Monte Carlo noise $\mathcal{E}_M$ and the discretization error $\mathcal{E}_D$. The Monte Carlo noise in replacing $\bar{\rho}(\mathbf{y}^{j},h)$ by $\rho_{M}(\mathbf{y}^j,h)$ is identical to the Monte Carlo noise in PIC, that is, $\mathcal{O}(\sqrt{\frac{\rho(\mathbf{y}^j)}{Nh^D}})$. The discretization error in the first step has two sources: the Taylor expansion in (\ref{eq:density_apcloud}) and the GFD approximation of the gradient in (\ref{eq:pseudoinverse}). The difference between the average of the $k$th order Taylor expansion of $\rho$ and $\rho$ itself is $\mathcal{O}(\rho(\mathbf{y}^j)h^{2\lfloor \frac{k}{2} \rfloor+2})$, where we obtain an additional order for even $k$ due to the symmetry of the kernel $\Phi$. Because the error of $l$th order derivative is $\mathcal{O}(\rho(\mathbf{y}^j)h^{k-l+1})$, and the coefficient for $l$th order derivative $a_{d_1,d_2,\cdots,d_l}(\mathbf{y}^j)=\mathcal{O}(h^{l})$, the discretization error given by the GFD derivative approximation is $\mathcal{O}(\rho(\mathbf{y}^j)h^{k+1})\geq \mathcal{O}(\rho(\mathbf{y}^j)h^{2\lfloor \frac{k}{2} \rfloor+2})$. Thus the total discretization error in step 1 is $\mathcal{O}(\rho(\mathbf{y}^j)h^{k+1})$.

Generalized finite-difference Poisson solver has the same accuracy with its estimation of $\Delta \phi$, i.e., $\mathcal{O}(\rho(\mathbf{y})h^{k-1})$. However, for the second order GFD we observe a supraconvergence when the GFD stencil is well-balanced due to the error cancellation similar to that in the standard five point finite difference stencil. 
In this case, the error for the solution and its gradient of the GFD Poisson solver are both $\mathcal{O}(f(\mathbf{y})h^2)$, as observed in our numerical experiments. The interpolation from $\mathcal{C}$ to $\mathcal{P}$, based on the $k$th order Taylor expansion of $\phi$ where the derivatives are given by GFD, is $(k+1)$th order accurate for $\phi$ and $k$th order accurate for $\mathbf{E}$.

\subsection{Refinement strategy for AP-Cloud}

Generally, when $k$th order GFD is used in the charge assignment scheme, the Poisson solver, and the differentiation and interpolation routines, 
the total error for both $\phi$ and $\mathbf{E}$ are
\begin{equation}
\mathcal{E}=\mathcal{O}(\mbox{E}[|\mathcal{E}_M|])+\mathcal{O}(\mathcal{E}_D)
=\mathcal{O}\left(\sqrt{\frac{\rho(\mathbf{y})}{Nh^D}}+\rho(\mathbf{y})h^{k-1}\right),
\end{equation}
where the $\mathcal{O}(\rho(\mathbf{y})h^{k-1})$ leading term in discretization error is from GFD Poisson solver.
To minimize the error, the optimal mesh size is
\begin{equation}\label{gfdoptimal}
h=\mathcal{O}\left(\frac{1}{N\rho(\mathbf{y})}\right)^{\frac{1}{2k+D-2}},
\end{equation}
and the minimized error is 
\begin{equation}\label{merror}
\mathcal{E}_{total}=\rho(\mathbf{y})^{\frac{k+D-2}{2k+D-2}}N^{-\frac{k-1}{2k+D-2}}.
\end{equation}

For the second order GFD in particular, we have better error bound due to the symmetry of interpolation kernel and stencil and supraconvergence
\begin{equation}
\mathcal{E}=\mbox{E}[|\mathcal{E}_M|]+\mathcal{E}_D
=\mathcal{O}\left(\sqrt{\frac{\rho(\mathbf{y})}{Nh^D}}+\rho(\mathbf{y})h^2\right),
\end{equation}
and the optimal mesh size is the same as in (\ref{eq:optimal})
\begin{equation}
h=\mathcal{O}\left({N\rho(\mathbf{y})}\right)^{-\frac{1}{4+D}}.
\end{equation}

\section{Implementation}
\label{sec:implementation}
We use a $2^D$-tree data structure to store particles, and select some of its cell centres as computational nodes. The $2^D$-tree data structure is a tree data structure in a $D$-dimensional space in which each cell has at most $2^D$ children. Quadtree and octree are standard terms in 2D and 3D spaces, respectively. 

The algorithm in~\cite{octree} is used in the $2^D$-tree construction. The first step is to sort the particles  by their Morton key, so particles in the same cell are contiguous in the sorted array. Then leaf cells are constructed by an array traversal, during which we record the number of particles and the index of the first particle in each cell. At last the interior cells are constructed in a depth decreasing order by a traversal of cells of the deeper level. The overall time complexity, dominated by the Morton key sorting, is $\mathcal{O}(N\log N)$, where $N$ is the total number of particles. This parallel $2^D$-tree construction algorithm, together with parallel linear solver, enables efficient parallel implementation of AP-Cloud.

Because all computational nodes are cell centres of a $2^D$-tree, their distribution will be similar to an AMR-PIC mesh. This is a result of our implementation method and not an internal property of AP-Cloud.

\subsection{Error balance criterion}

The optimal interparticle distance $h$ given in (\ref{gfdoptimal}) depends on the charge density $\rho$. In most applications, we do not know $\rho$ in advance; otherwise, we do not need the charge assignment scheme to estimate it. We use a Monte Carlo method to obtain a rough estimation of $\rho$:
\begin{equation}\label{rough}
\rho(\mathbf{y})=\frac{N'}{NV},
\end{equation}
where $V$ is the volume of a neighborhood of $\mathbf{y}$, and $N'$ is the number of particles in the neighborhood.

If the neighborhood is the box with the edge length $h$ centred at $\mathbf{y}$, $V=h^D$, we substitute (\ref{rough}) into (\ref{gfdoptimal}) and obtain 
\begin{equation}\label{criterion}
h=\mathcal{O}({N'}^{-\frac{1}{2k-2}}).
\end{equation}

\subsection{2:1 mesh balance}

If the charge density undergoes rapid changes, as is typical for certain applications (such as particle accelerators and cosmology), the optimal $h$ (\ref{criterion}) also changes rapidly.
This causes two potential problems. First, when the difference between levels of refinement on two sides of a cell is too large, that cell cannot find a balanced GFD stencil. If no particle in the coarse side is chosen to be its neighbor, the information on that side is missing. If we force the algorithm to choose a particle on the coarse side as a neighbor, the truncation error from this particle is much larger than that from the others. 
Second, in some cases, there are almost no particles in the region near the boundary. 
In order to enforce the boundary condition, interior nodes need to use far 
away nodes located on the boundary as their neighbors. 

To avoid these two problems, we enforce a 2:1 mesh balance. The 2:1 mesh balance requires that the difference between the levels of refinement of two neighbors is at most one. Because the mesh size changes smoothly, both imbalanced GFD stencils and empty regions are avoided. 

\subsection{Algorithm to select computational nodes and search GFD neighbors}
For clarity, we will focus on the selection of nodes in 3D in this subsection. The selection of nodes in 2D is similar and easier. An octree cell and the centre of a cell will be used interchangeably in this subsection.

We say an octree cell $\mathbf{z}$ is a neighbor of another octree cell $\mathbf{y}$ in a set of octree cells $S$, if 
\begin{itemize} 
\item $\mathbf{y}, \mathbf{z}\in S$;
\item $Cell(\mathbf{y})\cap Cell(\mathbf{z})$ is a face.
\item $level(\mathbf{z})\leq level(\mathbf{y})$;
\item No ancestor of $\mathbf{z}$ satisfies the previous three conditions.
\end{itemize}

The neighbors defined here are different from the neighbors in~\cite{octree} or the neighbors in GFD stencil. Generally, for any cell $\mathbf{y}$ in any set of octree cells $S$, it has at most 6 neighbors, each corresponding to one of its 6 faces. It is possible for a cell to have less than 6 neighbors. For example, root cell has no neighbor in any $S$ in a non-periodic region.

During the selection of nodes, we will keep a queue of octree cells, $Q$, and a list $L(\mathbf{y})$ for each cell $\mathbf{y}$ in $Q$ containing its neighbors in $Q$, which we call neighbor list. The basic operation is to open a cell $\mathbf{y}\in Q$:
\begin{enumerate}
\item Mark $\mathbf{y}$ as non-node candidate.
\item Add all children cells of $\mathbf{y}$ at the end of $Q$, mark them as node candidates.
\item Initialize the neighbor lists of the new added cells. Some of the neighbors are their siblings, while others are the neighbors of $\mathbf{y}$ or children of the neighbors of $\mathbf{y}$.
\item Update the neighbor lists of the descendent of the neighbors of $\mathbf{y}$.
\end{enumerate}

The algorithm for the selection of nodes is as follows.

\begin{enumerate}
\item Initialize a queue $Q$, which contains only the root cell.
\item Traverse $Q$. For each cell $\mathbf{y}$ in $Q$, test if it satisfies
\begin{equation}\label{errorbalance_implementation}
h<c{N'}^{-\frac{1}{2k-2}},
\end{equation}
where $c$ is a tuning parameter, $h$ is the diameter of the subtree, $n$ is the number of particles in the subtree. If the condition is not satisfied, open $\mathbf{y}$. Let $l$ be the deepest level in $Q$ at the end of this traverse.
\item Traverse $Q$. For each leaf cell $\mathbf{y}$ at level $l$, check if the neighbors of $\mathbf{y}$ satisfy the 2:1 mesh balance. Open each neighbor $\mathbf{z}$ that does not satisfy 2:1 mesh balance.
\item If $l>0$, let $l\leftarrow l-1$, then repeat step 3. If $l=0$, output all node candidates as computational nodes. If non-periodic boundary condition is used, add additional nodes on the boundary.
\end{enumerate}

Given the number of particles and the order of GFD is fixed, the tuning parameter $c$ in (\ref{errorbalance_implementation}) determines the number of nodes. Ideally, $c$ can be computed from the constant in the proportional relationship in (\ref{gfdoptimal}), which in turn depends on the order of GFD, the kernel function, and the relative magnitude of $\rho$ and its gradients. However, in most applications, the relative magnitude of $\rho$ and its gradients is unknown, so we try different values of $c$ and compare their results to estimate its optimal value in numerical tests.

Checking error balance criterion and 2:1 mesh balance takes only constant number of operations per cell. Except the part to update the neighbor lists, each open operation takes constant number of elementary operations as well, which can be charged on the 8 new added cells, so the time complexity is $\mathcal{O}(|Q|)$, where $|Q|$ is the number of cells in $Q$ in the end of the selection. To analyze the complexity to update the neighbor lists, we note that each time we update the neighbor of a cell $\mathbf{z}$ in $Q$, the level of its neighbor increases. Since the level of its neighbor is bounded by the height of the octree $l_{max}$, the total running time to update neighbors of all particles is $\mathcal{O}(|Q|l_{max})$. Because each interior cell of $Q$ has 8 children, $Q$ is a complete octree, we have $|Q|<8/7n$, where $n$ is the number of not opened cells, i.e., the number of computational nodes. In conclusion, the complexity to select nodes is $\mathcal{O}(|Q|l_{max})$.

After selecting the nodes, the neighbor list can be used to search GFD neighbors. If $\mathbf{z}$ is a neighbor of node $\mathbf{y}$, and $\mathbf{z}$ is not a node itself, then the 4 children of $\mathbf{z}$ that share a face with $\mathbf{y}$ must be nodes because of 2:1 mesh balance. The nodes among the neighbors of $\mathbf{y}$ and the children of these neighbors that share a face with $\mathbf{y}$ are called 1-ring.  The union of the $k$-ring and the nodes among the neighbors of $k$-ring nodes and the children of these neighbors that share a face with $k$-ring nodes are called $(k+1)$-ring. GFD neighbors are chosen from 2-ring by the quadrant criterion if there are enough number of nodes in 2-ring. If there are not enough number of nodes in $k$-ring, we will try to select GFD neighbors from $(k+1)$-ring. In our simulations, 5-ring always contains enough neighbors. Because this neighbor searching algorithm only depends on the local nodes information and local data structures, we claim the complexity to search neighbors of a node is independent from the total number of nodes, and this algorithm takes $\mathcal{O}(n)$ time to find neighbors of all computational nodes.

One problem is related to the fact that we do not know in advance how deep the nodes are in the octree while we build the octree.
In other words, it is possible that during the algorithm to select nodes, a leaf cell in the octree needs to be opened. However, this is very unlikely to happen in practice, if we always use the maximum depth supported in the implementation. For example, in the 2D Gaussian beam with halo test in Section 5, the order of GFD is 2, the region is $[-1,1]\times[-1,1]$, and the minimum tuning parameter $c$ used in simulations is 0.01. If we use a 64 bits Morton key, the maximum depth of the octree is 21, the cell size of the leaf cell is $(1-(-1))\times 2^{-21}=2^{-20}$. According to (\ref{errorbalance_implementation}), $0.01^2\times 2^{20\times 2}=1.10\times 10^8$ particles must be in the same leaf cell in order to open it, which is more than the total number of particles in the whole domain.

\section{Numerical results}

The AP-Cloud method approximately minimizes the error. However, it does incur additional
cost for the construction and search of the octree data structure. In addition, the 
fast Fourier transform can no longer be used for solving the resulting linear system,
so we must replace it with a sparse linear solver. Therefore, the practical advantage
of AP-Cloud is by no means obvious. In this section, we present some numerical results
for problems with non-uniform distribution, and demonstrate the advantages of AP-Cloud
in terms of both accuracy and efficiency compared to PIC. We also discuss potential advantages over AMR-PIC.

\subsection{2D Gaussian beam with halo}

We have performed verification of the Adaptive Particle-in-Cloud method using examples of highly non-uniform distributions of particles typical for accelerator beams with halos. 
In such problems,  a high-intensity, small-sigma particle beam is surrounded by a larger radius halo containing from 3 to 6 orders of magnitude smaller number of particles compared
to the main beam. As accurate modeling of realistic accelerator beam and halo distributions is unnecessary for the numerical verification, we represent the system by axially symmetric
Gaussian distributions. This also allows us to obtain a benchmark solution.

Consider  the following 2D electrostatic problem 
\begin{equation}\label{ep}
\Delta \phi=\rho,
\end{equation} 
where charge density $\rho$ is given by two overlapping Gaussian distribution:
\begin{equation}
\rho(\mathbf{x})=a_1\left[\exp\left(-\frac{|\mathbf{x}|^2}{\tau_1^2}\right)+a_2\exp\left(-\frac{|\mathbf{x}|^2}{\tau_2^2}\right)\right],
\label{GaussianDist}
\end{equation} 
in the domain $\Omega=[-1,1]\times [-1,1]$. We use the following values for the coefficients: the radius of the main beam $\tau_1=0.02$, the halo
intensity  $a_2=10^{-5}$, and the width of the halo $\tau_2=0.3$. Coefficient $a_1=396.1$ is a normalization parameter to ensure 
$\int_\Omega\rho(\mathbf{x}) d\mathbf{x}=1$. The model is consistent (in terms of the order of magnitude for the beam versus halo ratio) with real particles beams  in accelerators. 

While the AP-Cloud method is independent of the geometric shape of the computational domain, we solve the problem in a square domain to enable the comparison with the 
traditional PIC method. The benchmark solution is obtained in the following way. The problem is embedded in a larger domain, a radius 2 disk, using the same 
 charge density function $\rho$  and the homogeneous Dirichlet boundary condition.  A solution, obtained by a highly refined 1D solver in cylindrically symmetric coordinates,
 is considered as the benchmark solution. The Dirichlet  boundary condition for the two-dimensional problem is computed by interpolating the 1D solution at the location of the 2D boundary.
 This boundary condition function is then used for both the second-order AP-Cloud and PIC methods.

In our numerical simulations, CIC scheme~(\ref{cloud_in_cell}) is used in charge assignment and interpolation in PIC method. Theoretically, there are more accurate schemes available, such as triangular shaped cloud with reshaping step. However, these higher order schemes are very computationally intensive, and are not able to give better result than CIC with the same CPU time in our numerical tests. The order of accuracy of AP-Cloud method does not depend on the particular kernel function $\Phi$, so we choose the nearest grid point scheme for its simplicity, that is, $\Phi$  in~(\ref{eq:density_apcloud}) is set to be the characteristic function of the corresponding octree cell. 

\begin{figure}[H]
\begin{center}
\includegraphics[width=8cm]{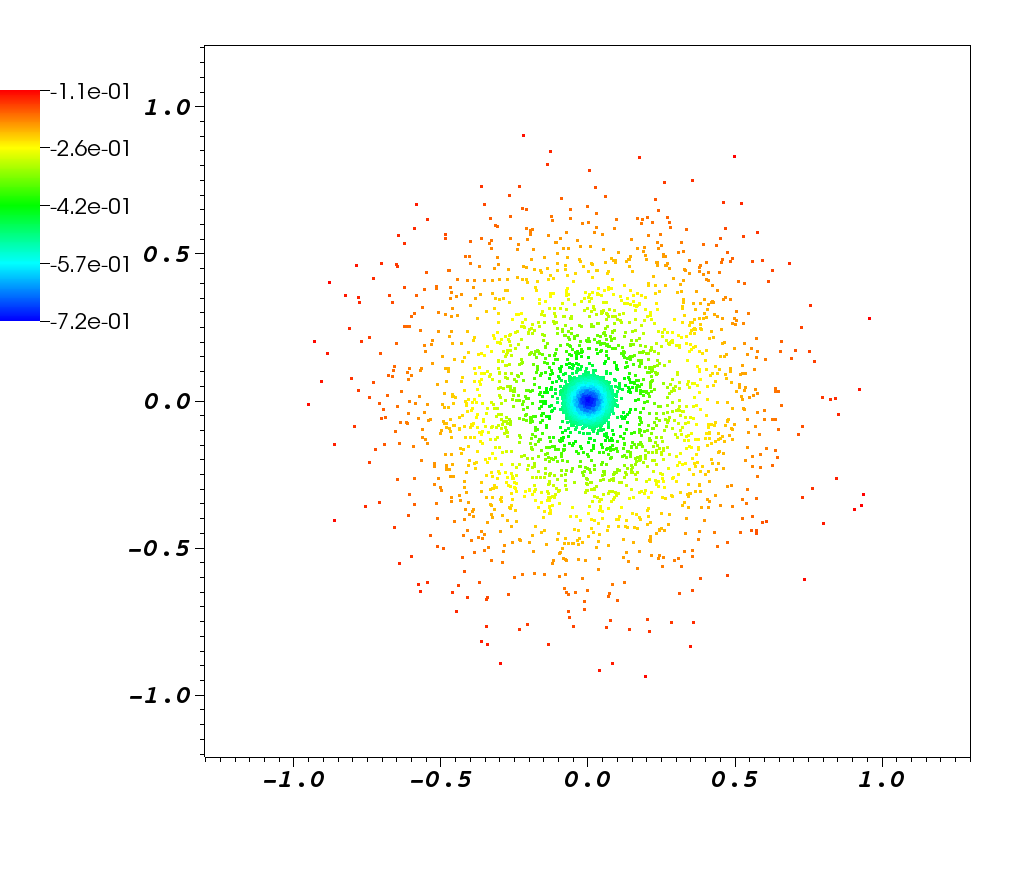}
\caption{Distribution of $10^6$ particles and the electric potential.}
\label{2D_physics_particles}
\end{center}
\end{figure}

\begin{figure}[H]
\begin{center}
\includegraphics[width=8cm]{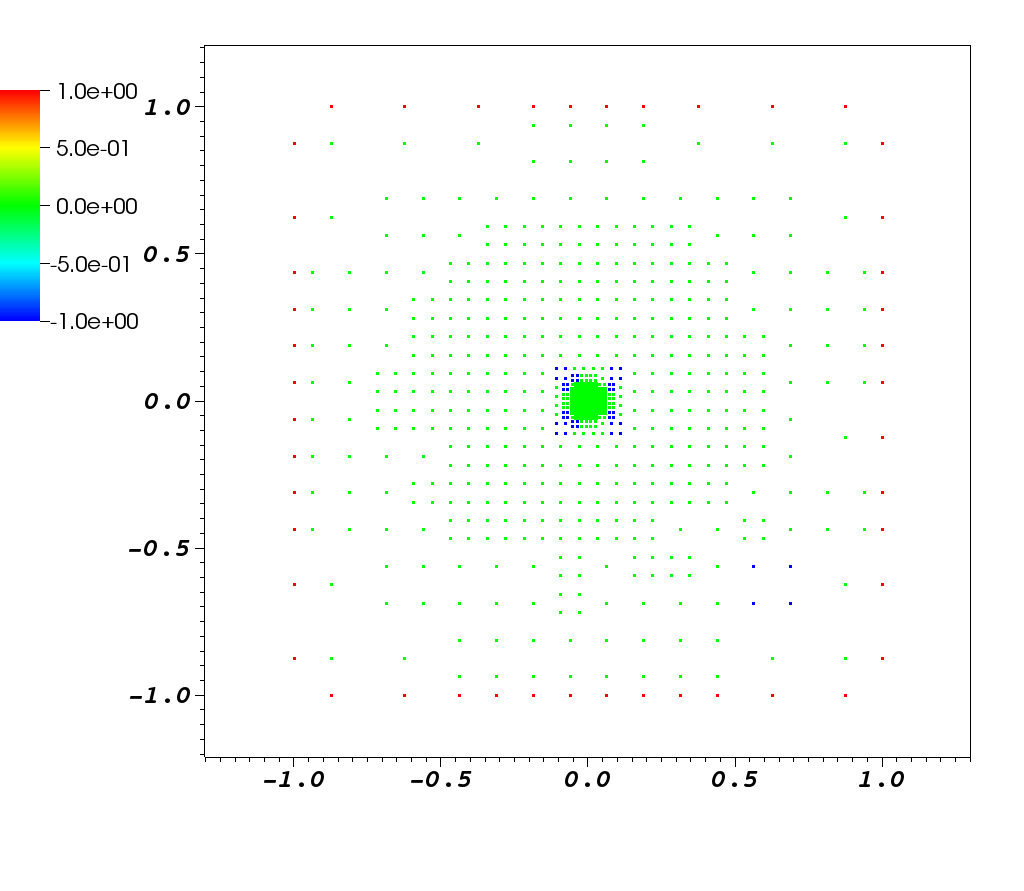}
\caption{Distribution of nodes.  Green nodes are given by error balance criterion, blue nodes are given by 2:1 mesh balance, and red nodes are on the boundary.}
\label{2D_comput_particles}
\end{center}
\end{figure}

\begin{figure}[H]
\begin{center}
\includegraphics[width=8cm]{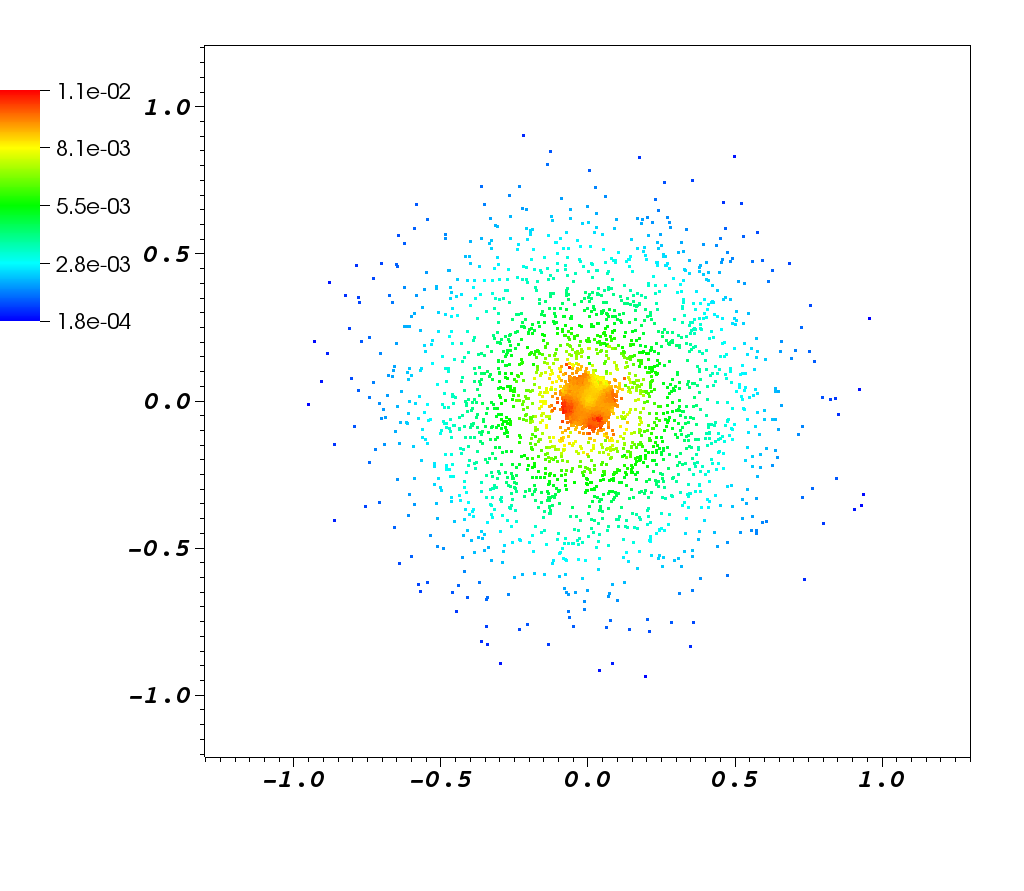}
\caption{Distribution of error of potential $\phi$ computed by AP-Cloud. ($\Vert\phi\Vert_\infty=0.7226.$)}
\label{2D_ap_cloud_error}
\end{center}
\end{figure}

\begin{figure}[H]
\begin{center}
\includegraphics[width=8cm]{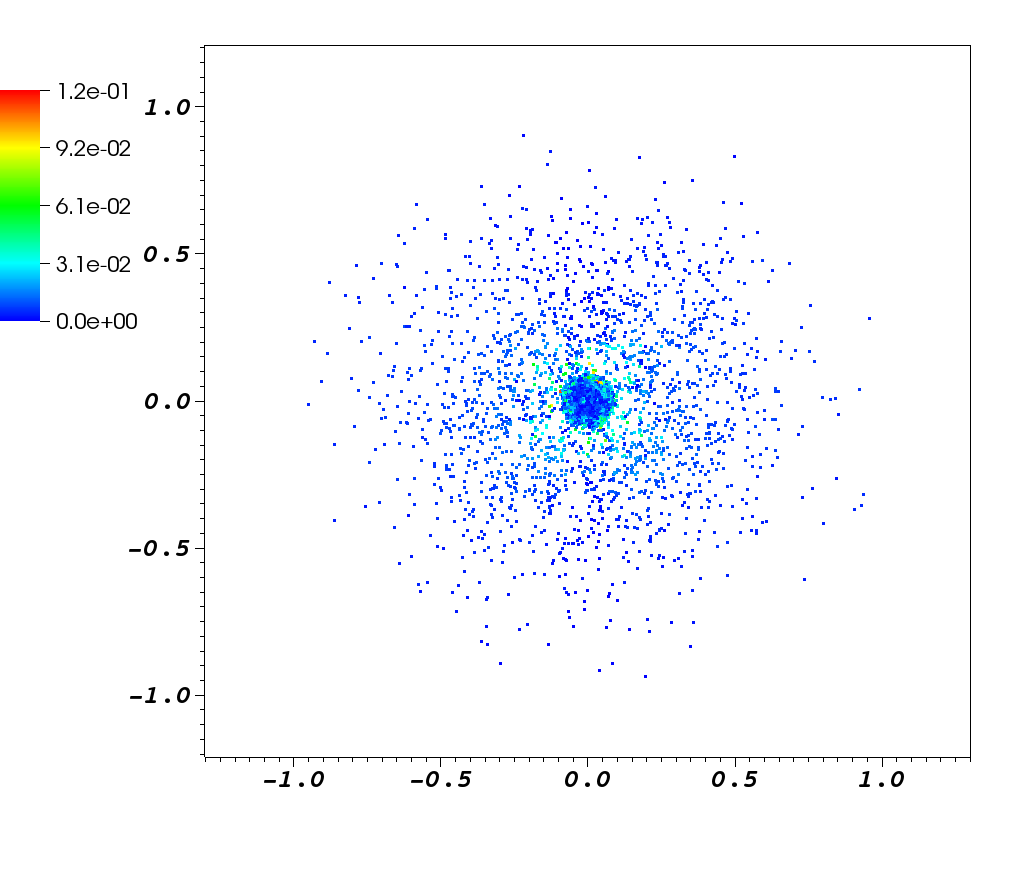}
\caption{Distribution of error of gradient of potential $\phi_x$ computed by AP-Cloud. ($\Vert\phi_x\Vert_\infty=3.581.$)}
\label{2D_ap_cloud_grad}
\end{center}
\end{figure}

\begin{table}[H]
\begin{center}
\begin{tabular}{|c|c|c|c|}
\hline
n&Running time&Error of $\phi$&Error of $\phi_x$
\\ \hline
121	&0.341&	0.118& 	2.05 \\ \hline
441	&0.327&	0.0648& 	1.87 \\ \hline
1681	&0.353&	0.0347 &	1.41 \\ \hline
6561	&0.490&	0.0139 &	0.674 \\ \hline
25921&	1.60&	0.00371& 	0.214 \\ \hline
\end{tabular}
\caption{CPU time and accuracy of traditional PIC with $10^6$ particles.}
\label{2D_PIC}
\end{center}
\end{table}

\begin{table}[H]
\begin{center}
\caption{CPU time and accuracy of AP-Cloud with 2:1 mesh balance with $10^6$ particles.}
\label{2D_ap_cloud}
\begin{tabular}{|c|c|c|c|}
\hline
n&Running time&Error of $\phi$&Error of $\phi_x$
\\ \hline
256	&0.929&	0.00289& 	0.0515 \\ \hline
428	&0.933&	0.0183 &	0.0218 \\ \hline
1156	&0.970&	0.00886 &	0.00927 \\ \hline
3652	&1.14&	0.00365  &	0.00750  \\ \hline
7559	&1.45&	0.000233& 	0.00725 \\ \hline
19077&2.80&	0.000145& 	0.00724 \\ \hline
\end{tabular}
\end{center}
\end{table}

\begin{figure}[H]
\begin{center}
\includegraphics[width=8cm]{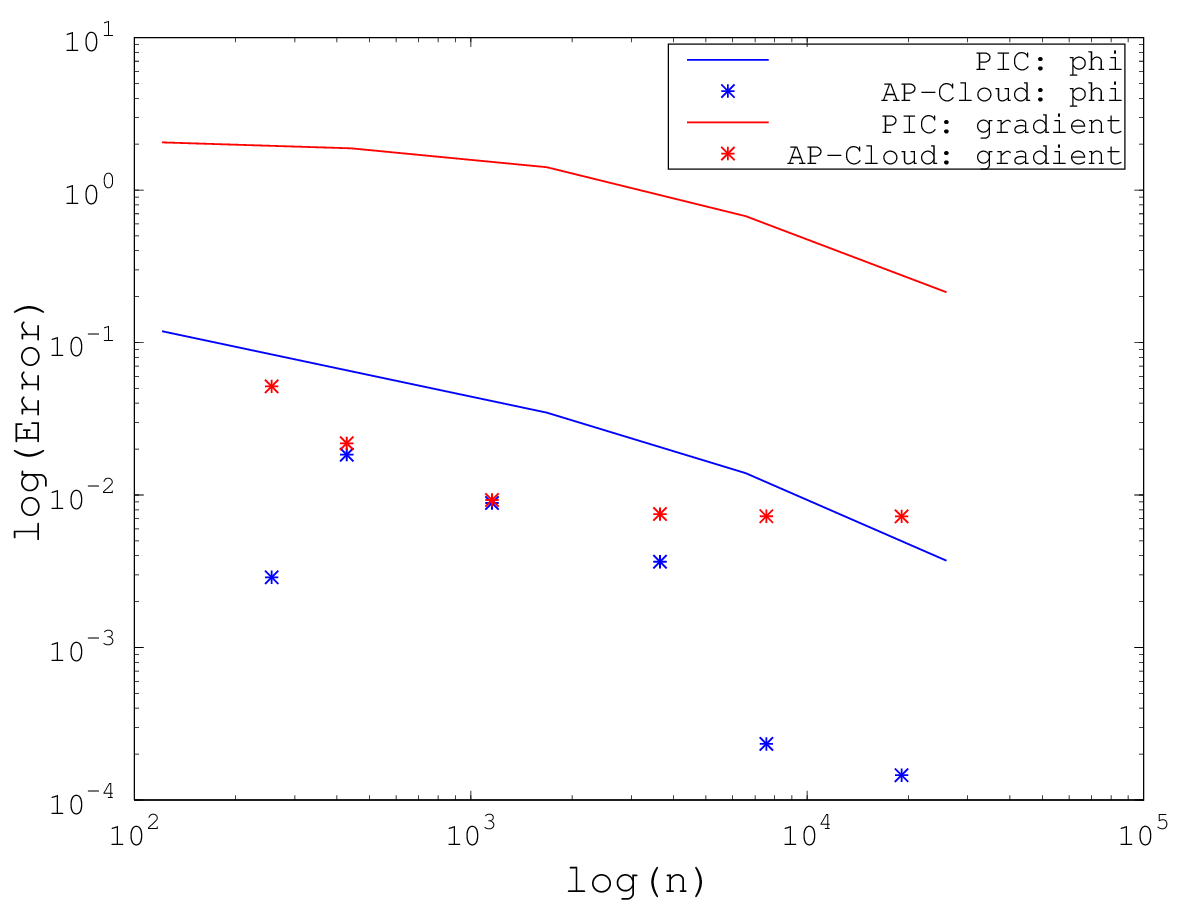}
\caption{Dependence of error on the number of nodes for PIC and AP-Cloud.}
\label{log_plot_ap_cloud}
\end{center}
\end{figure}

\begin{figure}[H]
\begin{center}
\includegraphics[width=8cm]{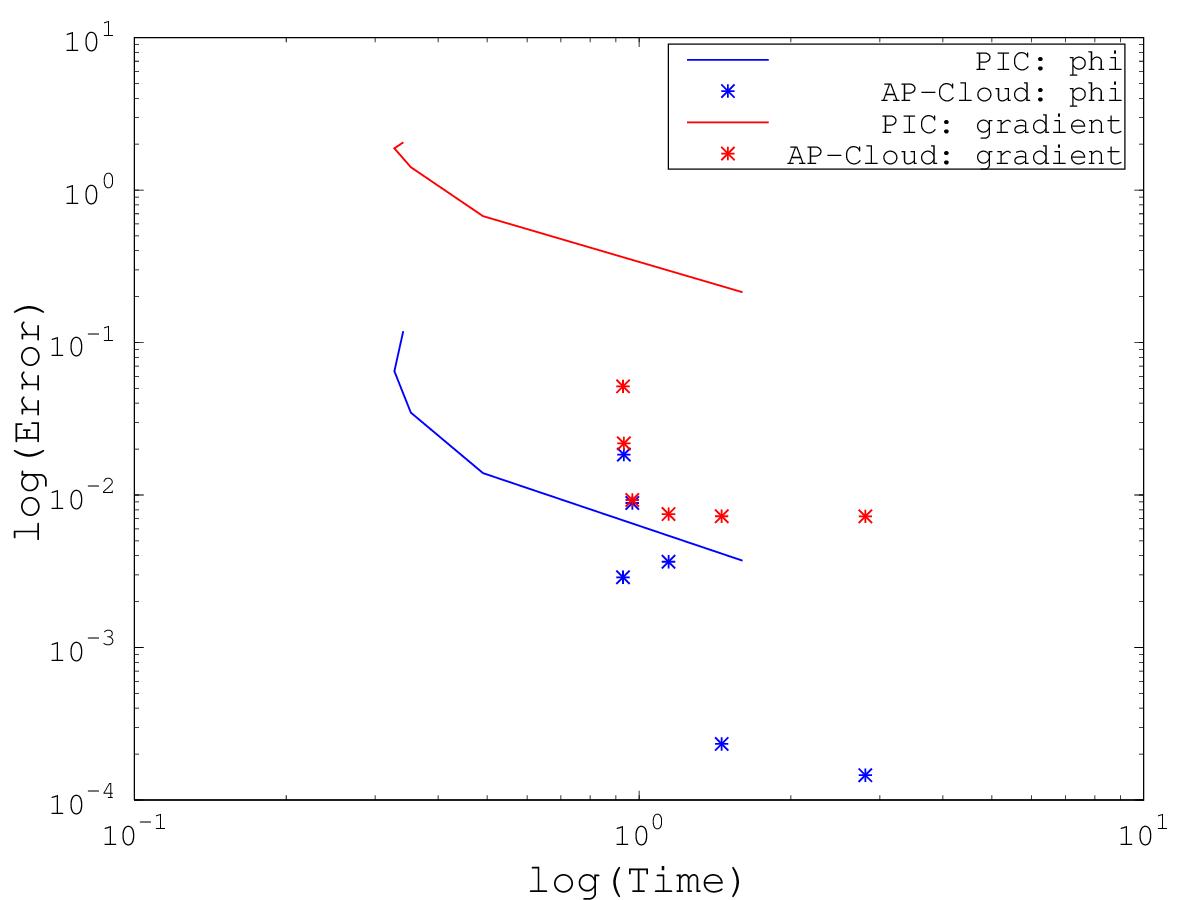}
\caption{Dependence of error on CPU time for PIC and AP-Cloud.}
\label{log_plot_pic}
\end{center}
\end{figure}

 Figures \ref{2D_physics_particles} - \ref{2D_ap_cloud_grad} show the distribution of particles coloured according to solution values, nodes, 
 and distributions of errors of the potential and its gradient. $L_2$ errors are used in Tables \ref{2D_PIC} and \ref{2D_ap_cloud} and Figures \ref{log_plot_ap_cloud} and \ref{log_plot_pic} show that the estimation of the potential and its
gradient given by AP-Cloud is much more accurate compared to the PIC estimation. For example, the gradient error computed by AP-Cloud with 256 nodes is only about one fourth of the error of PIC with 19077 nodes. Although AP-Cloud is computationally more intensive for the same number of nodes due to the construction of a 
quadtree and solving an additional linear system for $\rho$, its accuracy under the same 
running time is still significantly better.

\begin{table}[H]\label{tab:runningtime}
\begin{center}
\caption{Breakdown of running times of AP-Cloud.}
\begin{tabular}{|c|c|c|c|c|c|}
\hline
N&$10^4$&$10^5$&$10^6$&$10^6$&$10^6$\\ \hline
n&828&1900&1888& 4420& 11190\\ \hline
Build quadtree&  7.67e-03&5.48e-02&  4.85e-01&   4.85e-01&   4.89e-01\\ \hline
Search nodes& 3.11e-04&6.58e-04&   6.46e-04&   1.45e-03&   3.70e-03\\ \hline
Build linear systems&  1.11e-02&   2.44e-02&  2.39e-02&   5.62e-02&   1.41e-01\\ \hline
Solve linear system for $\rho$ &  1.64e-01&1.68e-01&  1.68e-01&   1.95e-01&   2.76e-01\\ \hline
Solve linear system for $\phi$ &  1.81e-01&1.95e-01&  1.95e-01&   2.99e-01&   7.93e-01\\ \hline
Find interpolation coefficient &    1.41e-02 &  3.20e-02 &  3.15e-02 &  7.46e-02&   1.88e-01\\ \hline
Interpolate &8.18e-04&9.11e-03&  1.40e-01 &  1.40e-01 &  1.45e-01 \\ \hline
Total running time &3.82e-01&   4.89e-01  &   1.05e+00  & 1.26e+00 &  2.05e+00 \\ \hline
\end{tabular}
\end{center}
\end{table}

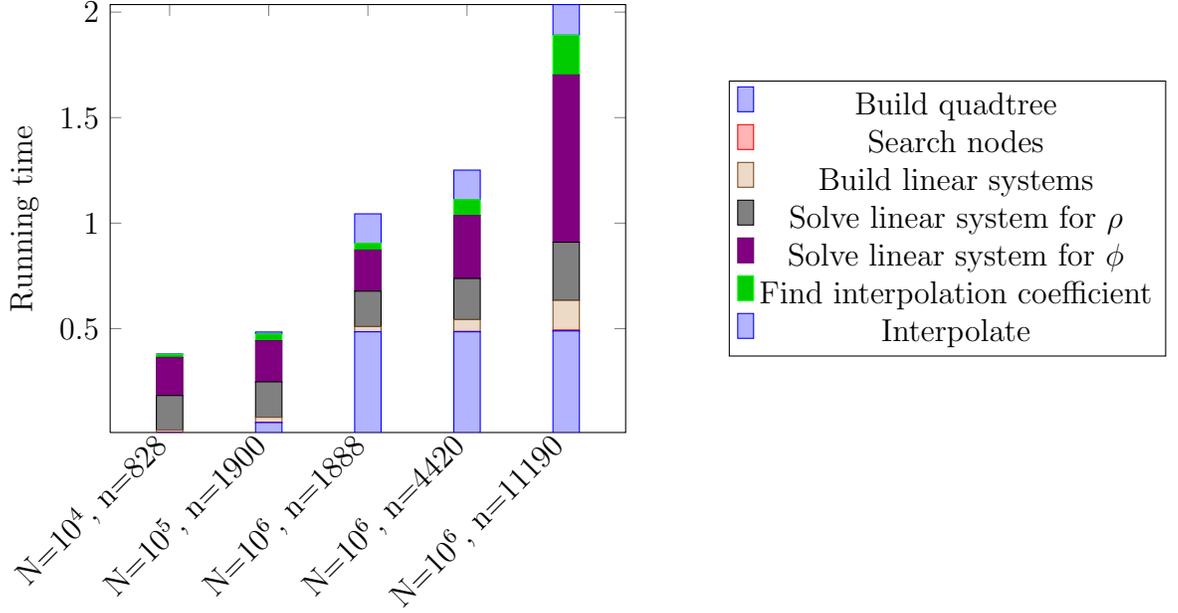
\begin{figure}
\begin{tikzpicture} 
 \begin{axis}[ 
ybar stacked, 
enlargelimits=0.15, 
enlarge y limits=false, 
legend style={at={(1.2,0.5)}, 
anchor=west,legend columns=1}, 
ylabel={Running time}, 
symbolic x coords={{N=$10^4$, n=828}, {N=$10^5$, n=1900}, {N=$10^6$, n=1888}, {N=$10^6$, n=4420}, {N=$10^6$, n=11190}}, 
xtick=data, 
x tick label style={rotate=45,anchor=east}, 
] 
 \addplot+[ybar] plot coordinates {({N=$10^4$, n=828},7.67e-03) ({N=$10^5$, n=1900},5.48e-02) ({N=$10^6$, n=1888},4.85e-01) ({N=$10^6$, n=4420},4.85e-01) ({N=$10^6$, n=11190},4.89e-01)};
 \addplot+[ybar] plot coordinates {({N=$10^4$, n=828},3.11e-04) ({N=$10^5$, n=1900},6.58e-04) ({N=$10^6$, n=1888},6.46e-04) ({N=$10^6$, n=4420},1.45e-03) ({N=$10^6$, n=11190},3.70e-03)};
 \addplot+[ybar] plot coordinates {({N=$10^4$, n=828},1.11e-02) ({N=$10^5$, n=1900},2.44e-02) ({N=$10^6$, n=1888},2.39e-02) ({N=$10^6$, n=4420},5.62e-02) ({N=$10^6$, n=11190},1.41e-01)};
 \addplot+[ybar] plot coordinates {({N=$10^4$, n=828},1.64e-01) ({N=$10^5$, n=1900},1.68e-01) ({N=$10^6$, n=1888},1.68e-01) ({N=$10^6$, n=4420},1.95e-01) ({N=$10^6$, n=11190},2.76e-01)};
  \addplot+[ybar] plot coordinates {({N=$10^4$, n=828},1.81e-01) ({N=$10^5$, n=1900},1.95e-01) ({N=$10^6$, n=1888},1.95e-01) ({N=$10^6$, n=4420},2.99e-01) ({N=$10^6$, n=11190},7.93e-01)};
  \addplot+[ybar] plot coordinates {({N=$10^4$, n=828},1.41e-02) ({N=$10^5$, n=1900},3.20e-02) ({N=$10^6$, n=1888},3.15e-02) ({N=$10^6$, n=4420},7.46e-02) ({N=$10^6$, n=11190},1.88e-01)};
  \addplot+[ybar] plot coordinates {({N=$10^4$, n=828},8.18e-04) ({N=$10^5$, n=1900},9.11e-03) ({N=$10^6$, n=1888},1.40e-01) ({N=$10^6$, n=4420},1.40e-01) ({N=$10^6$, n=11190},1.45e-01)};

 \legend{Build quadtree, Search nodes, Build linear systems, Solve linear system for $\rho$, Solve linear system for $\phi$,  Find interpolation coefficient, Interpolate } 
 \end{axis} 
\end{tikzpicture}
\caption{Stacked plot of running times of AP-Cloud.}
\label{stack}
\end{figure}

From theoretical complexity analysis in Section 4 and experimental results in Table \ref{tab:runningtime} and Figure \ref{stack}, the steps of AP-Cloud can be divided into 3 main groups:
\begin{itemize}
\item Quadtree construction and interpolation, which time complexity is $\mathcal{O}(N\log N)$. The running time for this group is $\mathcal{O}(N\log N)$, which dominates when $N\gg n$.
\item Searching for nodes, building linear systems and finding interpolation coefficients. The running time for this group is $\mathcal{O}(nl_{max})$, which is small compared to the running time of other two groups.
\item Solving the linear system for $\rho$ and $\phi$. CPU time depends on both the linear solver and $n$, and dominates for small ratios of $N/n$.
\end{itemize}

In this test, we did not obtain the second order convergence due to the Monte Carlo noise. Table \ref{convergence_order} shows the result of another test, where $\bar{\rho}(\mathbf{y}^{j},h)$ in (\ref{eq:density}) is given by the integral of the exact density function instead of Monte Carlo integration. The convergence is second order for both the potential and gradient, as expected.

\begin{table}[H]
\begin{center}
\caption{Convergence of AP-Cloud without Monte Carlo noise.}
\label{convergence_order}
\begin{tabular}{|c|c|c|c|c|}
\hline
n&Error of $\phi$&Error of $\phi_x$&Order of $\phi$&Order of $\phi_x$
\\ \hline
240& 0.00181& 0.00123& -& -\\ \hline
863 &0.000564& 0.000315&    1.82  & 2.13\\ \hline
3336 &0.000150& 7.62e-05&   1.95 &  2.09\\ \hline
13043 &3.80e-05& 2.006e-05 &  2.01&   1.96\\ \hline
\end{tabular}
\end{center}
\end{table}

\subsection{3D Gaussian beam with halo}

\begin{figure}[H]
\begin{center}
\includegraphics[width=8cm]{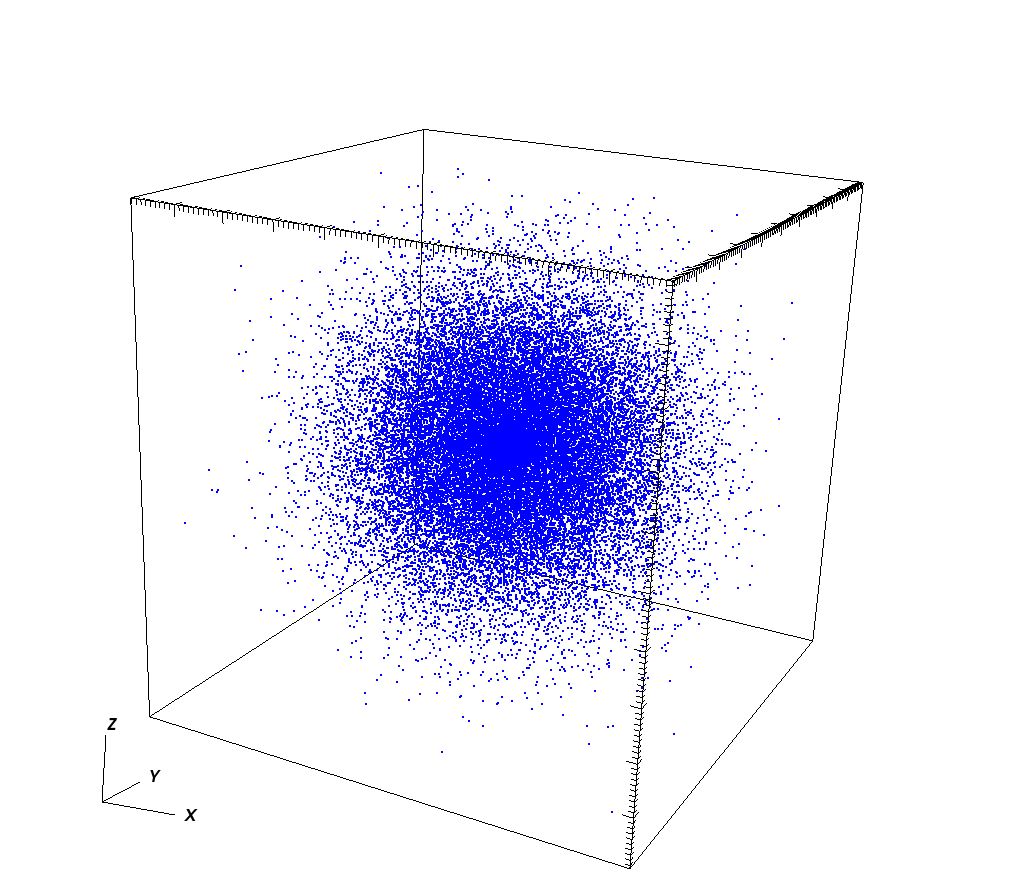}
\caption{Distribution of $10^6$ particles.}
\label{3D_particles}
\end{center}
\end{figure}

\begin{figure}[H]
\begin{center}
\includegraphics[width=8cm]{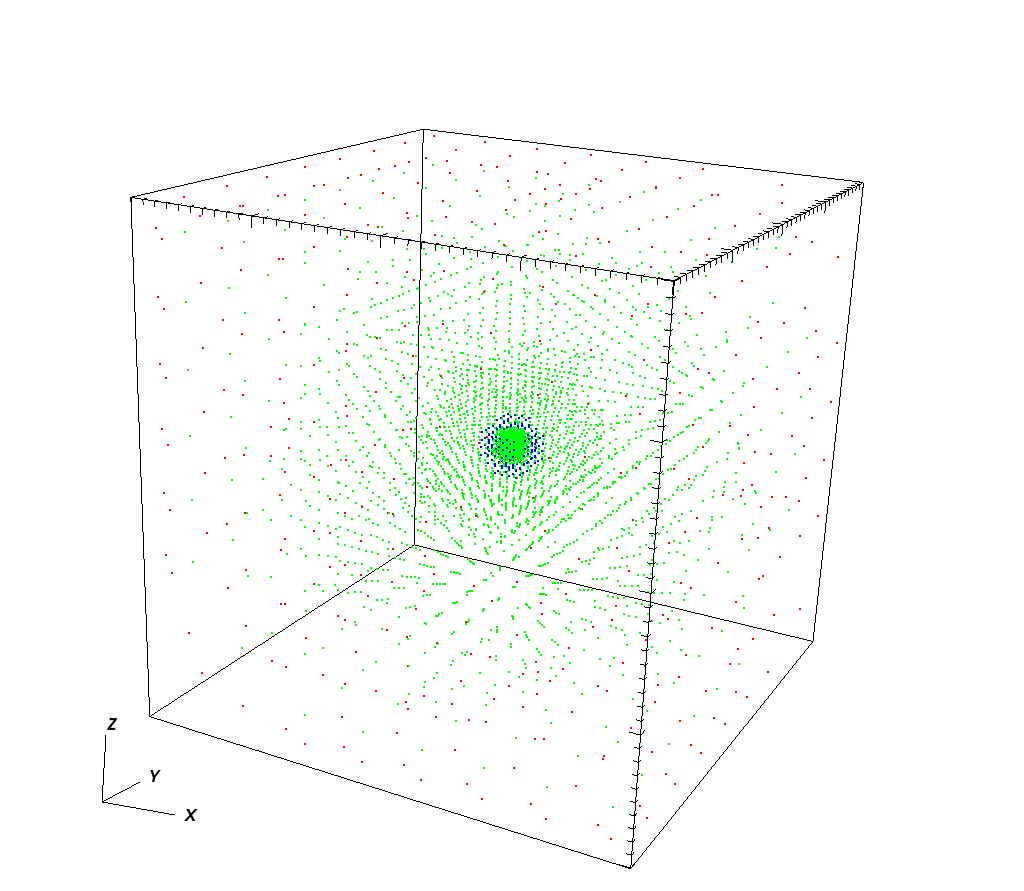}
\caption{Distribution of 4067 nodes. Green nodes are given by error balance criterion, blue nodes are given by 2:1 mesh balance, and red nodes are on the boundary.}
\label{3D_comput_particles}
\end{center}
\end{figure}

In this Section, we investigate the accuracy of the AP-Cloud method in
3D. To enable comparison with a simple benchmark solution, we study a
spherically symmetric extension of the beam-with-halo problem. Despite
the loss of physics relevance, it is a useful problem that tests
adaptive capabilities of the method. Consider the Poisson equation
with the charge density $\rho$ given by two overlapping Gaussian
distributions (\ref{GaussianDist}) in the domain
$\Omega=[-1,1]^3$. The radius of center beam is $\tau_1=0.02$, the
strength of the halo is $a_2=10^{-5}$, and the width of the halo is
$\tau_1=0.3$. The coefficient $a_1=7677$ provides the normalization
$\int_\Omega\rho(\mathbf{x}) d\mathbf{x}=1$. The benchmark solution
and the boundary condition function were obtained similarly to the 2D
case. The distribution of $10^6$ particles is shown in Figure
\ref{3D_particles}.  The AP-Cloud computation, performed using 4067
nodes (Figure \ref{3D_comput_particles}), gives a solution with the
normalized norm of $\phi$ on particles is 0.0352
($\Vert\phi\Vert_\infty=3.038$), and the normalized norm of $\phi_x$
on particles is 0.578 ($\Vert\phi_x\Vert_\infty=41.18$).

\begin{figure}[H]
\begin{center}
\includegraphics[width=8cm]{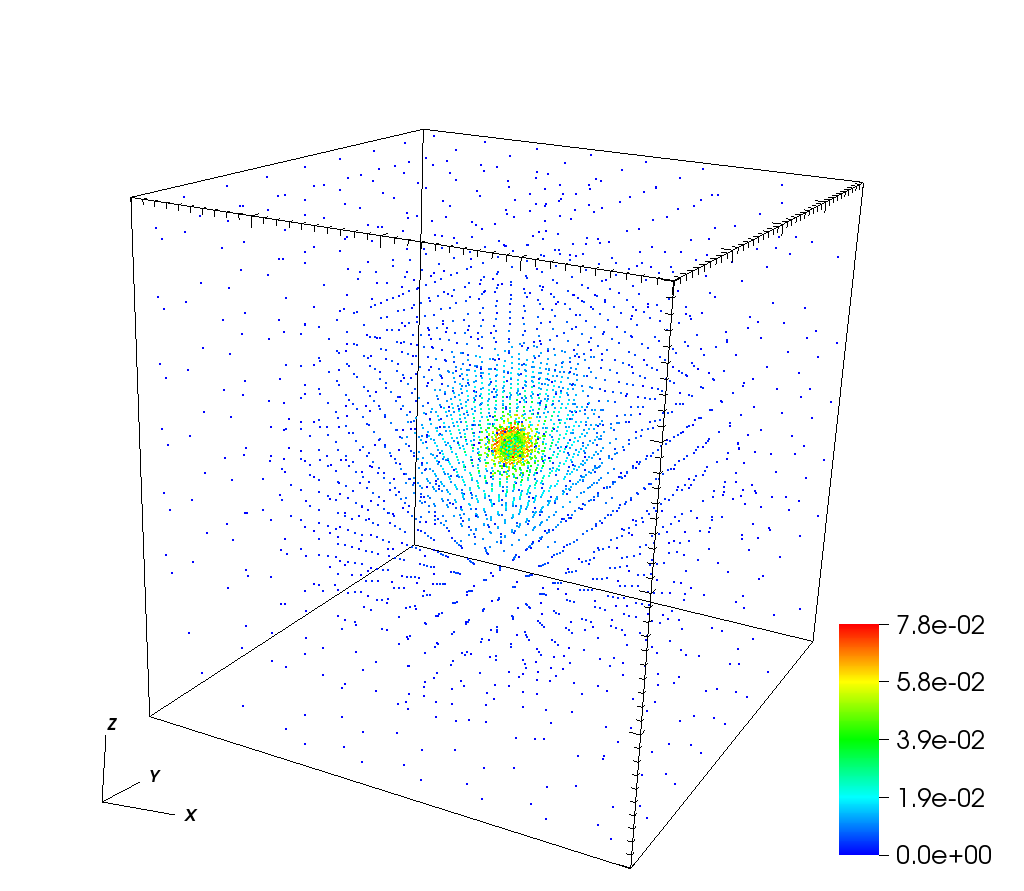}
\caption{Distribution of error of potential $\phi$ by AP-Cloud ($\Vert\phi\Vert_\infty=3.038$).}
\end{center}
\end{figure}

\begin{figure}[H]
\begin{center}
\includegraphics[width=8cm]{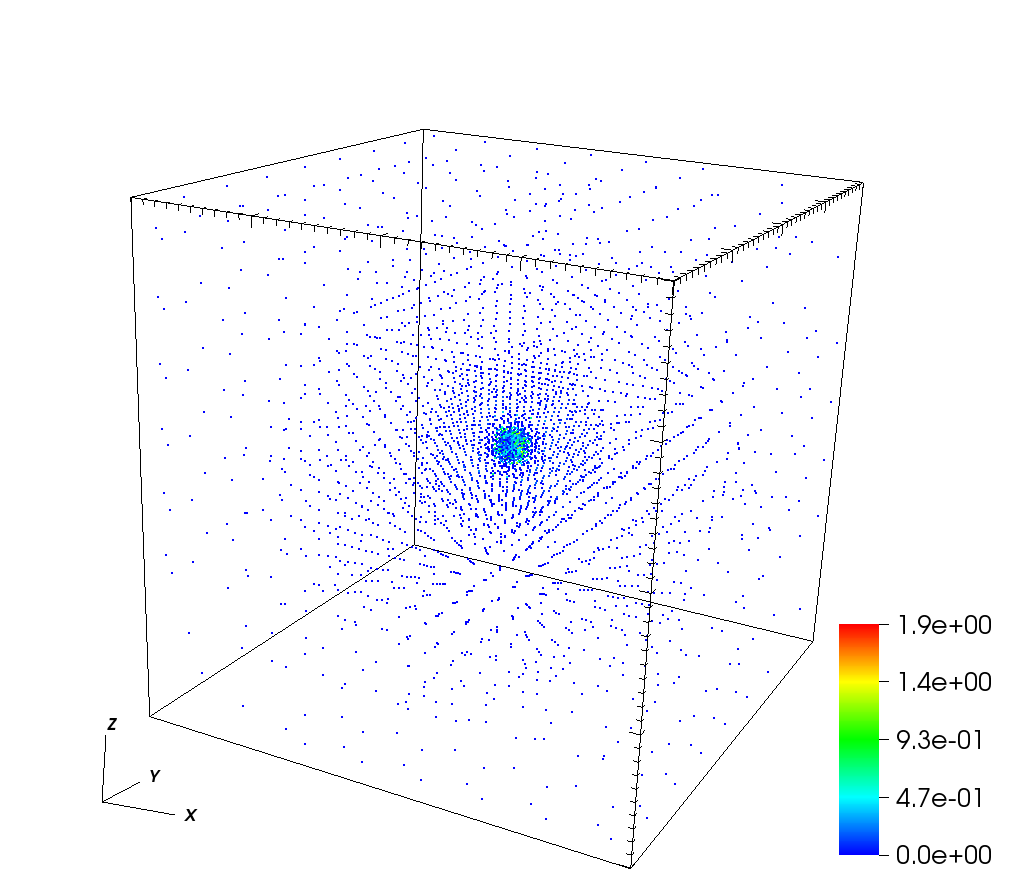}
\caption{Distribution of error of gradient of potential $\phi_x$ by AP-Cloud ($\Vert\phi_x\Vert_\infty=41.18$).}
\end{center}
\end{figure}


\begin{table}[H]
\begin{center}
\begin{tabular}{|c|c|c|c|}
\hline
n&Running time&Error of $\phi$&Error of $\phi_x$
\\ \hline
8000	&0.443&	1.40& 	19.6 \\ \hline
64000	&1.48&	0.726& 	16.1 \\ \hline
512000&24.1&	0.219 &	8.60 \\ \hline
4096000&361&	0.0606 &	2.92 \\ \hline\end{tabular}
\caption{CPU time and accuracy  of traditional PIC with $10^6$ particles in 3D.}
\label{3D_PIC}
\end{center}
\end{table}

\begin{table}[H]
\begin{center}
\begin{tabular}{|c|c|c|c|}
\hline
n&Running time&Error of $\phi$&Error of $\phi_x$
\\ \hline
1546	&0.921&	0.0402& 	1.17 \\ \hline
4067	&1.14&	0.0352 &	0.578 \\ \hline
13687	&2.10&	0.0183 &	0.329 \\ \hline
59349&7.22&	0.00443  &	0.244  \\ \hline
\end{tabular}
\caption{CPU time and accuracy of AP-Cloud with 2:1 mesh balance with $10^6$ particles in 3D.}
\label{3D_ap_cloud}
\end{center}
\end{table}

\begin{figure}[H]
\begin{center}
\includegraphics[width=8cm]{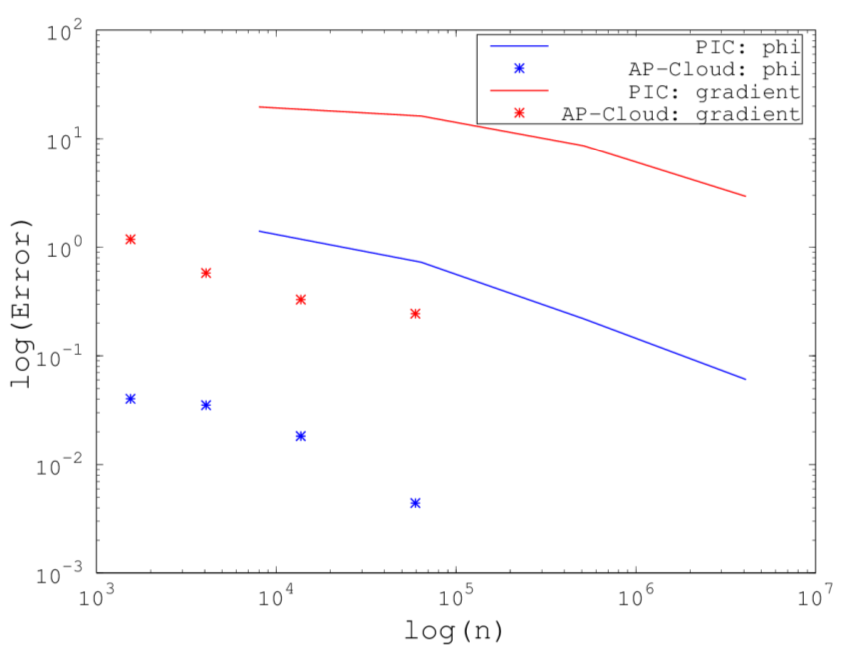}
\caption{Dependence of error on the number of nodes for PIC and AP-Cloud in 3D.}
\label{log_plot_3d_n}
\end{center}
\end{figure}

\begin{figure}[H]
\begin{center}
\includegraphics[width=8cm]{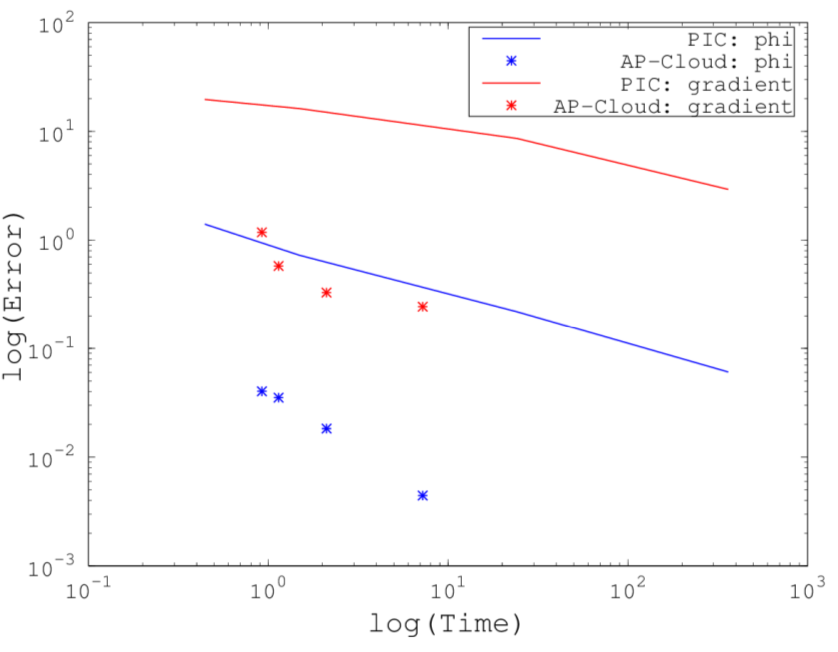}
\caption{Dependence of error on CPU time for PIC and AP-Cloud in 3D.}
\label{log_plot_3d_t}
\end{center}
\end{figure}

Results presented in Tables \ref{3D_PIC} and \ref{3D_ap_cloud} and
Figures \ref{log_plot_3d_n} and \ref{log_plot_3d_t} show that the
advantage of AP-Cloud is more evident in 3D problem. With only 1546
nodes and less than 1 second running time, AP-Cloud has more accurate
result than PIC with 4096000 cells and more than 361 seconds running
time.


\subsection{Test for self-force effect with single particle}

As mentioned in the introduction, Vlasov-Poisson problems with highly
non-uniform distributions of matter can be solved using the adaptive
mesh refinement technique for PIC \cite{VayCol1,VayCol2}. However, it
is well known that AMR-PIC introduces significant artifacts in the
form of artificial image particles across boundaries between coarse
and fine meshes. These images introduce spurious forces that may
potentially alter the particle motion to an unacceptable level
\cite{VayCol1,VayCol2}.  Methods for the mitigation of the spurious
forces have been designed in \cite{ColNor10}. The traditional PIC on a
uniform mesh is free of such artifacts.

The convergence of Adaptive Particle-in-Cloud solutions to benchmark
solutions, discussed in the previous Section, already indicates the
absence of artifacts. To further verify that AP-Cloud is free of
artificial forces present in the original AMR-PIC, we have performed
an additional test similar to the one in \cite{VayCol1}, which
involved the motion of a single particle across the coarse and fine
mesh interface. For AP-Cloud, we studied the motion of a single test
particle represented by a moving cloud of nodes with refined distances
towards the test particle. The test particle contained a smooth,
sharp, Gaussian-type charge distribution to satisfy the requirements
of the GFD method.

\begin{figure}[H]
\begin{center}
\includegraphics[width=0.49\textwidth]{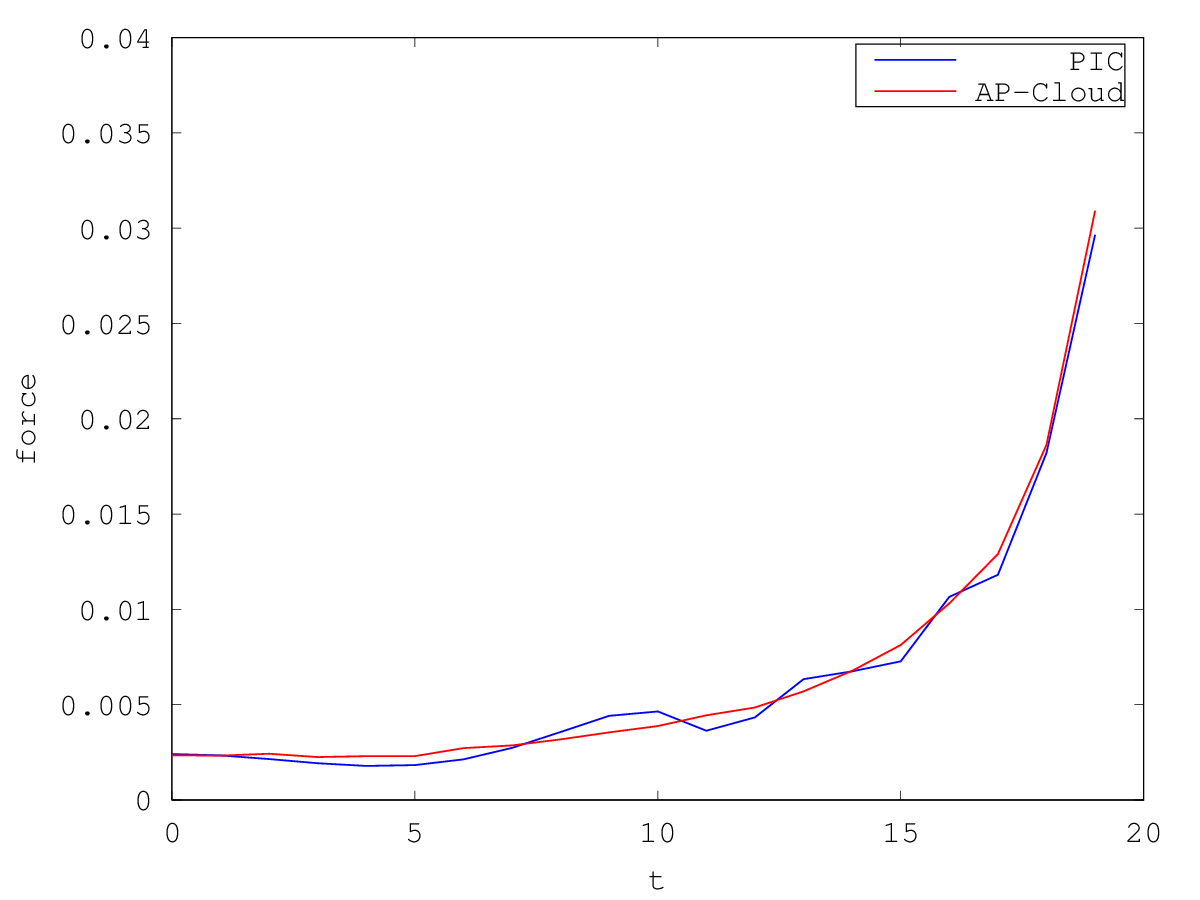}
\includegraphics[width=0.49\textwidth]{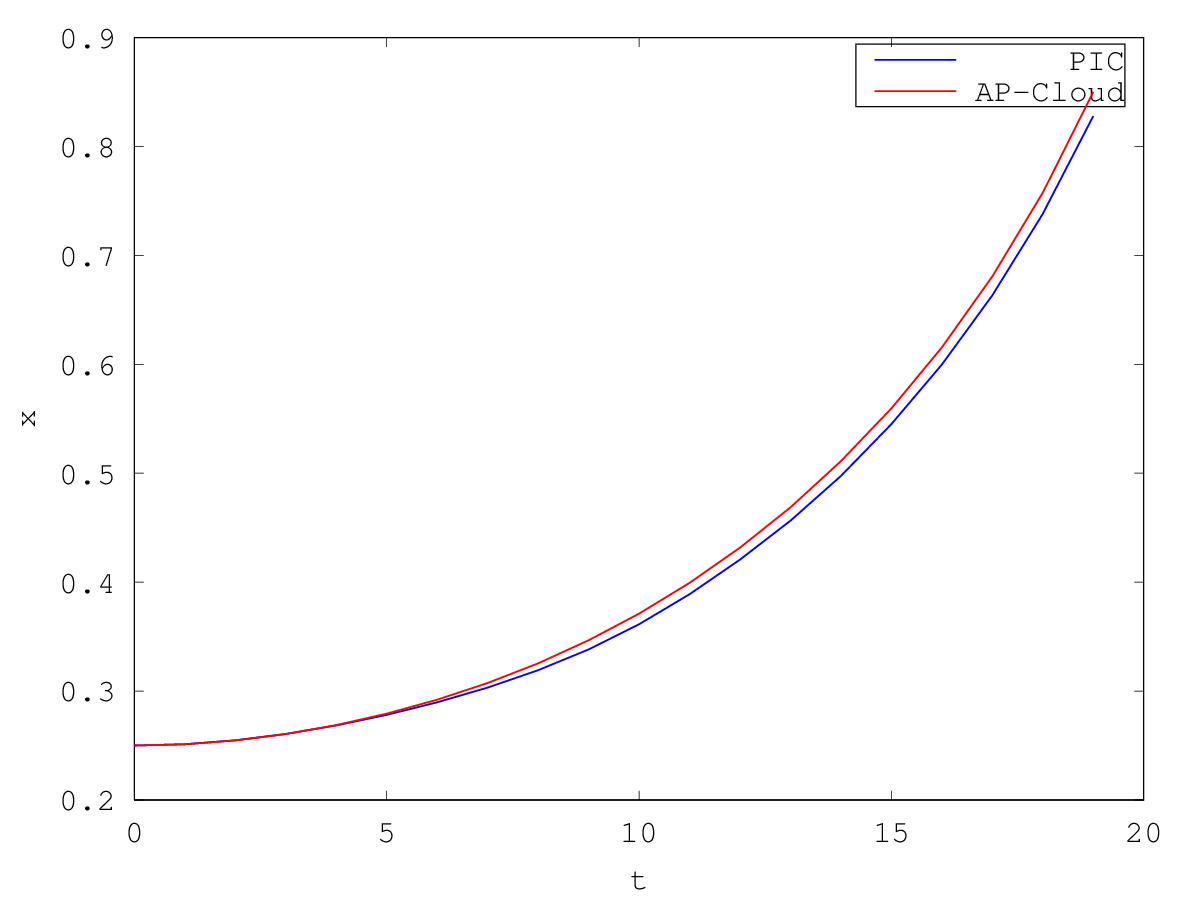}
\caption{Forces (left) and motion (right) of a single test particle
  obtained with PIC and AP-Cloud methods demonstrating the absence of
  artifacts in the AP-Cloud method.}
\label{fig:single_p_trajectory}
\end{center}
\end{figure}

The forces and motion of a single test particle obtained with PIC and
AP-Cloud methods are shown in Figure~\ref{fig:single_p_trajectory}. We
observe that the electric forces computed by the AP-Cloud method are
more accurate and smoother compared to even the traditional PIC.  But
the oscillatory deviation of forces in PIC from the correct direction
does not cause accumulation of the total error due to conservative
properties of PIC.  The trajectories of the particle obtained by both
methods are close.  The test provides an additional assurance that
artificial images are not present in the AP-Cloud method.

\section{Summary and Conclusions}

We have developed an Adaptive Particle-in-Cloud (AP-Cloud) method that
replaces the Cartesian grid in the traditional PIC with adaptive
computational nodes.  Adaptive particle placement balances the errors
of the differential operator discretization and the source computation
(analogous to the error of the Monte Carlo integration) to minimize
the total error.

AP-Cloud uses GFD based on weighted least squares (WLS) approximations
on a stencil of irregularly placed nodes. The framework includes
interpolation, least squares approximation, and numerical
differentiation capable of high order convergence.

The adaptive nature of AP-Cloud gives it significant advantages over the
traditional PIC for non-uniform distributions of particles and complex
boundaries. It achieves significantly better accuracy in the gradient
of the potential compared to the traditional PIC for the problem of
particle beam with halo. The method is independent of the geometric
shape of the computational domain, and can achieve highly accurate
solutions in geometrically complex domains. The optimal mesh size based on error-balance criterion gives AP-Cloud a potential 
advantage over AMR-PIC in terms of accuracy, and specially designed tests
showed that the AP-Cloud method is free of artificial images and
spurious forces typical for the original AMR-PIC without special
mitigation techniques. Another advantage of AP-Cloud over AMR-PIC is 
the ease of implementation, as AP-Cloud does not require special remapping routines between different meshes. Our future work will focus on 
higher convergence rates of the method, performance optimization, 
parallel implementation using hybrid technologies, as well as applications to practical problems with non-uniform distribution of matter. 
A direct comparison of AP-Cloud with AMR-PIC in terms of accuracy and  efficiency will also be addressed in the future work.

\section*{Acknowledgement}
This work was supported in part by the U.S. Department of Energy, Contract No. DE-AC02-98CH10886.

\appendix
\section{Error analysis for PIC}\label{appen}
Assume the interpolation kernel $\Phi$ is symmetric, non-negative, bounded by 1,
integrable and normalized over a compact support. That is, 
$0\leq \Phi(\mathbf{x})=\Phi(-\mathbf{x})\leq \Vert\Phi\Vert_\infty\leq 1$ for all 
$\mathbf{x}\in \mathbb{R}^D$,
$\Phi(\mathbf{x})=0$ for $\Vert\mathbf{x}\Vert>r$ for some $r$, and
$\int_{\mathbb{R}^D}\Phi(\mathbf{x})\mathbf{dx}=1$. 
These properties hold for all commonly used charge assignment schemes, including the 
nearest grid point, cloud-in-cell, and triangular shaped cloud schemes without reshaping step. In the following 
subsections, we will analyze the errors in the three steps of PIC, respectively.

\subsection{Error in Step 1}
To analyze the error in (\ref{eq:density}), let us first define an average quantity
at a grid point $\mathbf{y}^{j}$ as
$$\bar{\rho}(\mathbf{y}^{j},h):=\frac{1}{h^{D}}\int_{\Omega}\rho(\mathbf{x})\Phi\left(\frac{\mathbf{x}-\mathbf{y}^j}{h}\right)\mathbf{dx}.$$
Then,
$$\vert\tilde{\rho}(\mathbf{y}^{j},\mathcal{P},h)-\rho(\mathbf{y}^{j})\vert\leq\underbrace{\vert\tilde{\rho}(\mathbf{y}^{j},\mathcal{P},h)-\bar{\rho}(\mathbf{y}^{j},h)\vert}_{\mathcal{E}_{M}}+\underbrace{\vert\bar{\rho}(\mathbf{y}^{j},h)-\rho(\mathbf{y}^{j})\vert}_{\mathcal{E}_{D_1}},$$
where $\mathcal{E}_{M}$ is analogous to the error in Monte Carlo integration of a continuous
function, and $\mathcal{E}_{D_1}$ is known as moment error or discretization error in step 1, which depends on the interpolation kernel $\Phi$.
\subsubsection{Monte Carlo noise}

To bound the error $\mathcal{E}_{M}$, note that the expected value of $\tilde\rho$ is
\begin{equation}
\mbox{E}[\tilde\rho(\mathbf{y}^j,\mathcal{P},h)]=\bar\rho(\mathbf{y}^j,h).
\end{equation}
Therefore,
$$\mbox{Var}(\tilde\rho(\mathbf{y}^j,\mathcal{P},h))
=\mbox{E}[|\mathcal{E}_M|^2].$$

By the definition of variance,
\begin{equation}
\mbox{Var}(\tilde\rho(\mathbf{y}^j,\mathcal{P},h))=
\frac{1}{N}\mbox{Var}(\tilde\rho(\mathbf{y}^j,\mathbf{p}^i,h)),
\end{equation}
where $$\tilde\rho(\mathbf{y}^j,\mathbf{p}^i,h)=
\frac{1}{h^{D}}\Phi\left(\frac{\mathbf{p}^i-\mathbf{y}^j}{h}\right)$$ 
is the estimation of $\rho$ from a single particle $\mathbf{p}^i$, and
\begin{equation}
\mbox{E}[\tilde\rho(\mathbf{y}^j,\mathbf{p}^i,h)]=\bar\rho(\mathbf{y}^j,h).
\end{equation}
Let 
$\tau(\mathbf{y}^j,\mathbf{p}^i,h)=
\tilde\rho(\mathbf{y}^j,\mathbf{p}^i,h)-\bar\rho(\mathbf{y}^j,h).$
Then,
\begin{eqnarray}
\mbox{Var}(\tilde{\rho}(\mathbf{y}^{j},\mathbf{p}^{i},h))&=&\mbox{E}[\tau^{2}(\mathbf{y}^j,\mathbf{p}^i,h)]\\&=&\int_{\Omega}\tau^{2}(\mathbf{y}^j,\mathbf{x},h)\rho(\mathbf{x})\mathbf{dx}.
\end{eqnarray}
Consider the small neighborhood about $\mathbf{y}^j$ of radius $rh$, i.e., $\mathcal{A}^j=\{\mathbf{x}: \Vert\mathbf{x}-\mathbf{y}^j\Vert<rh\}$, so that 
$\tilde\rho(\mathbf{y}^j,\mathbf{x},h)=0$ in $\Omega\backslash\mathcal{A}^j$. Then,
$$
\mbox{Var}(\tilde{\rho}(\mathbf{y}^{j},\mathbf{p}^{i},h))=\underbrace{\int_{\mathcal{A}^{j}}\left(\tau(\mathbf{y}^{j},\mathbf{x},h)\right)^{2}\rho(\mathbf{x})\mathbf{dx}}_{I_{1}}+\underbrace{\int_{\Omega\backslash\mathcal{A}^{j}}\left(\tau(\mathbf{y}^{j},\mathbf{x},h)\right)^{2}\rho(\mathbf{x})\mathbf{dx}}_{I_{2}}.$$
Note that
\begin{equation}
\lim_{h\rightarrow 0}\bar\rho(\mathbf{y}^j,h)=\rho(\mathbf{y}^j)\int_{\mathbb{R}^D}\Phi(\mathbf{x})\mathbf{dx}
=\rho(\mathbf{y}^j)
\end{equation}
and
\begin{equation}
\lim_{h\rightarrow 0}\frac{1}{h^D}\int_{A^j}\rho(\mathbf{x})\mathbf{dx}=r^D V_D\rho(\mathbf{y}^j),
\end{equation}
where $V_D$ is the volume of the unit ball in $D$ dimensions.
Assume $h$ is small enough, so that $\bar\rho(\mathbf{y}^j,h)<2\rho(\mathbf{y}^j)$ and $\frac{1}{h^D}\int_{{A}^j}\rho(\mathbf{x})\mathbf{dx}<2r^D V_D\rho(\mathbf{y}^j)$. 
Then,
\begin{equation}
I_1\leq\int_{\mathcal{A}^{j}}\left(\frac{1}{h^{D}}\right)^{2}\rho(\mathbf{x})\mathbf{dx}<2\frac{r^{D}V_{D}}{h^{D}}\rho(\mathbf{y}^j)=\mathcal{O}(h^{-D})\rho(\mathbf{y}^j),
\end{equation}
and
\begin{equation}
I_2=\int_{\Omega\backslash\mathcal{A}^{j}}\bar{\rho}(\mathbf{y}^{j},h)^{2}
\rho(\mathbf{x})\mathbf{dx}<\bar{\rho}(\mathbf{y}^{j},h)^{2}<4\rho(\mathbf{y}^{j})^{2}=\mathcal{O}(1)\rho(\mathbf{y}^j).
\end{equation}
Therefore, \begin{eqnarray}
\mbox{Var}(\tilde\rho(\mathbf{y}^j,\mathcal{P},h))&=&\frac{1}{N}(I_1+I_2)\nonumber \\
&=&\frac{1}{N}\left(\rho(\mathbf{y}^j)\mathcal{O}\left(h^{-D}\right)+\rho(\mathbf{y}^j)^2\mathcal{O}(1)\right)\nonumber \\
&=&\frac{\rho(\mathbf{y}^j)}{N}\mathcal{O}\left(h^{-D}\right),
\end{eqnarray}
and 
\begin{equation}\label{montecarlo}
\mbox{E}[|\mathcal{E}_M|]=\sqrt{\mbox{E}[|\mathcal{E}_M|^2]-\mbox{Var}(|\mathcal{E}_M|)}\leq \sqrt{\mbox{E}[|\mathcal{E}_M|^2]}=\mathcal{O}\left(\sqrt{\frac{\rho(\mathbf{y}^j)}{Nh^D}}\right).
\end{equation}

Note that above bound for $\sqrt{\mbox{E}[|\mathcal{E}_M|^2]}$ is asymptotically tight. To see this, 
assume $4h^D\bar\rho(\mathbf{y}^j,h)<\Vert\Phi\Vert_\infty$, and let 
$A_1=\{\mathbf{x}: \Phi((\mathbf{x}-\mathbf{y}^j)/h)>\Vert\Phi\Vert_\infty/2\}$. Then,
$$I_{1}\geq\int_{A_{1}}\tau^2(\mathbf{y}^j,\mathbf{x},h)\rho(\mathbf{x})\mathbf{dx}>\int_{A_{1}}\left(\frac{\Vert\Phi\Vert_\infty}{4h^{D}}\right)^{2}\rho(\mathbf{x})\mathbf{dx}.$$
The volume of $A_1$ is $\mathcal{O}(h^D)$, so $I_{1}$ is no smaller than $\mathcal{O}(h^{-D})\rho(\mathbf{y}^j)$, and $\sqrt{\mbox{E}[|\mathcal{E}_M|^2]}$ is no smaller than $\mathcal{O}\left(\sqrt{\rho(\mathbf{y}^j)h^{-D}N^{-1}}\right)$.

\subsubsection{Moment error}

The moment error in charge assignment scheme, $\mathcal{E}_{D_1}$, can be viewed as the error in a numerical quadrature rule. Here we follow the approach taken in \cite{Cottet} and \cite{WangColella}. Let $B(\mathbf{y}^j,a)=\{\mathbf{x}:\Vert\mathbf{x}-\mathbf{y}^j\Vert<a\}$, $H(\mathbf{y}^j)$ be the Hessian matrix for $\rho$ at $\mathbf{y}^j$. The moment error is
\begin{eqnarray}\label{e1}
\mathcal{E}_{D_1}&=&\vert\bar{\rho}(\mathbf{y}^{j},h)-\rho(\mathbf{y}^{j})\vert\nonumber \\
&=&\frac{1}{h^D}\int_\Omega\left(\rho(\mathbf{x})-\rho(\mathbf{y}^j)\right)
\Phi\left(\frac{\mathbf{x}-\mathbf{y}^j}{h}\right)\mathbf{dx}   \nonumber \\
&=&\frac{1}{h^D}\int_{B(\mathbf{y}^j,rh)}\left((\mathbf{x-y}^j)^T\nabla\rho(\mathbf{y}^j)+(\mathbf{x-y}^j)^TH(\mathbf{y}^j)(\mathbf{x-y}^j) \right. \nonumber \\
&+&\left. \mathcal{O}(\Vert\mathbf{x-y}^j\Vert^3)\right)\Phi\left(\frac{\mathbf{x}-\mathbf{y}^j}{h}\right)\mathbf{dx}\nonumber \\
&=&\frac{1}{h^D}\int_{B(\mathbf{y}^j,rh)}(\mathbf{x-y}^j)^T\nabla\rho(\mathbf{y}^j)\Phi\left(\frac{\mathbf{x}-\mathbf{y}^j}{h}\right)\mathbf{dx} \nonumber \\
&+&\frac{1}{h^D}\int_{B(\mathbf{y}^j,rh)}(\mathbf{x-y}^j)^TH(\mathbf{y}^j)(\mathbf{x-y}^j)\Phi\left(\frac{\mathbf{x}-\mathbf{y}^j}{h}\right)\mathbf{dx}+\mathcal{O}(h^3)\nonumber \\
&=&0+\frac{1}{h^D}\int_{B(\mathbf{y}^j,rh)}\sum_{d=1}^D({x}_d-{y}_d^j)^2
\rho_{{y}_d{y}_d}(\mathbf{y}^j)
\Phi\left(\frac{\mathbf{x}-\mathbf{y}^j}{h}\right)\mathbf{dx}+\mathcal{O}(h^3)\nonumber \\
&=&h^2\sum_{d=1}^D\int_{B(\mathbf{y}^j,r)}({z}_d-{y}_d^j)^2\Phi\left(\mathbf{z}-\mathbf{y}^j\right)\mathbf{dz}\rho_{{y}_d{y}_d}(\mathbf{y}^j)+\mathcal{O}(h^3)\nonumber \\
&=&\sum_{d=1}^D\rho_{{y}_d{y}_d}(\mathbf{y}^j)\mathcal{O}(h^2),
\end{eqnarray}
where the fifth equal sign is due to the symmetry of interpolation kernel.
\subsection{Error in Step 2}
It is well known that the error in finite difference method for elliptic equation in step 2 is second order
\begin{equation}
\mathcal{E}_{D_2}=\sum_{d=1}^D\rho_{{y}_d{y}_d}(\mathbf{y}^j)\mathcal{O}(h^2).
\end{equation}

So the total discretization error in $\phi$ is
\begin{equation}
\mathcal{E_D}=\mathcal{E}_{D_1}+\mathcal{E}_{D_2}=\sum_{d=1}^D\rho_{{y}_d{y}_d}(\mathbf{y}^j)\mathcal{O}(h^2)=\rho(\mathbf{y}^j)\mathcal{O}(h^2),
\end{equation}
where the last equal sign is due to the assumption that $\rho$ and its
derivatives have comparable magnitude.
\subsection{Error in Step 3}
Step 3 contains a numerical differentiation and a linear
interpolation, both of which has second order accuracy. One problem is
that numerical differentiation is unstable. Generally, if $\phi$ is
already polluted with a $k$th order error, then its gradient given by
numerical differentiation has at most $(k-1)$th order accuracy,
because of the $h$ factor in denominator. For finite difference scheme
for elliptic equation, it has been proved that, nevertheless, if the
solution is smooth enough, then both the solution and its gradient are
second order convergent, even if the mesh is non-uniform, which is
called supraconvergence~\cite{sc}. More precisely, let $P_h$ be the
interpolation operator from grid function to piecewise linear
function, $R_h$ be the restriction operator from continuous function
to grid function, $u\in H^3(\Omega)$ be the exact solution of Poisson
equation, $u_h$ be the numerical solution, then we have
\begin{equation}
||P_h(R_hu-u_h)||_{H^1}\leq C h_{max}^2||u||_{H^3},
\end{equation} 
where $h_{max}$ is the maximal mesh size. Thus we claim the discretization error for $\mathbf{E}$ is also second order.

\end{document}